\documentclass[a4paper,12pt]{article}
\usepackage{a4,amssymb,amsmath,mathabx,hyperref,textgreek}
\usepackage{graphicx}   
\usepackage{color,float}
\usepackage[table]{xcolor}
\usepackage{graphics}
\usepackage{subfig}
\usepackage{tabularx}
\usepackage{blindtext}
\usepackage{authblk}
\graphicspath{{images/}}

\usepackage{epsfig}

\usepackage{tikz,pgfplots,filecontents,booktabs}
   \usetikzlibrary{shapes,shadows,arrows,matrix,trees,mindmap,backgrounds,automata,shapes.arrows,calc}
   \usetikzlibrary{patterns,decorations.markings,plotmarks}
\pgfplotsset{compat=1.8}
\usepackage[linesnumbered,boxed]{algorithm2e}
\definecolor{citecol}{rgb}{0.75,0,0}
\definecolor{urcol}{rgb}{0,0.5,0}
\definecolor{linkcol}{rgb}{0,0,0.75}
\definecolor{DarkRed}{rgb}{0.55,0.00,0.00}
\definecolor{grey}{gray}{0.95}

\hypersetup
{
bookmarksopen=true,
pdftitle="KMI",
pdfauthor="Anil kumar",
pdfmenubar=true, 
pdfhighlight=/O, 
colorlinks=false, 
pdffitwindow=true,     
pdfpagelayout=SinglePage, 
linkcolor=linkcol, 
citecolor=citecol, 
urlcolor=urcol 
}

\newtheorem{theorem}{Theorem}[section]
\newtheorem{remark}[theorem]{Remark}

\newtheorem{lemma}[theorem]{Lemma}

\newtheorem{definition}[theorem]{Definition}

\begin{document}

\setcounter{page}{1}

\title
{\Large \bf
Application of the Bell polynomials for the solution of some differential algebraic equations
}
\author[1,2,3,4]{Hari Mohan Srivastava \thanks{harimsri@math.uvic.ca}}
\author[5]{Giriraj Methi\thanks{girirajmethi@gmail.com}}
\author[5]{Anil Kumar\thanks{anil.digit@gmail.com }}
\author[6]{Mohammad Izadi\thanks{izadi@uk.ac.ir}}
\author[7]{Vishnu Narayan Mishra\thanks{vishnunarayanmishra@gmail.com} }
\author[8]{Brahim Benhammouda\thanks{bbenhammouda@hct.ac.ae}}

\affil[1]{\small{ Department of Mathematics and Statistics, University of Victoria,
	Victoria, British Columbia V8W 3R4,	Canada}}
\affil[2]{ Department of Medical Research, China Medical University Hospital, China Medical University, Taichung~40402, Taiwan, Republic of China}
\affil[3]{Department of Mathematics and Informatics, Azerbaijan University, 71 Jeyhun Hajibeyli Street, AZ1007~Baku, Azerbaijan}
\affil[4]{ Center for Converging Humanities, Kyung Hee University, 26 Kyungheedae-ro, Dongdaemun-gu, Seoul 02447, Republic of Korea}
\affil[5]{ Department of Mathematics \& Statistics, Manipal University Jaipur, Rajasthan, India}
\affil[6]{Department of Applied Mathematics, Faculty of Mathematics and Computer,
	Shahid Bahonar University of Kerman, Kerman, Iran}
\affil[7]{Department of Mathematics, Indira Gandhi National Tribal University, Madhya
	Pradesh,India}
\affil[8]{Higher Colleges of Technology, Abu Dhabi, United Arab Emirates}

\date{}
\maketitle
\begin{abstract}
	The differential transform method is used to find numerical approximation of solution to a class of certain nonlinear differential algebraic equations. The method is based on Taylor's theorem. Coefficients of the Taylor series are determined by constructing a recurrence relation. To deal with nonlinearity of the problems, the Fa\`{a} di Bruno's formula containing the partial ordinary Bell polynomials is applied within the differential transform to avoid computation of symbolic derivatives. The error estimation results are presented too. Four concrete problems are studied to show efficiency and reliability of the method. The obtained results are compared to other methods.

\end{abstract}
\begin{quotation}

\end{quotation}

\section{Introduction}

System of differential algebraic equations (DAEs) are combination of ordinary differential equations (ODEs) together with purely algebraic equations. Many researchers have investigated physical problems involving DAEs in electrical network \cite{Stri05}, modelling of constrained mechanical systems \cite{Ben22,Carmine20}, optimal control \cite{Brenan,Ghazwa21} and chemical processes problems \cite{Lio98}. \smallskip

\noindent DAEs with higher index $\left( > 1\right)$ are difficult to solve. Index reduction techniques can be used to convert them into lower index problems, but this is computationally expensive and sometimes changes properties of the solutions also. \smallskip

\noindent Several methods have been implemented to find solutions of DAEs such as Runge-Kutta method \cite{Linh18}, Adomian decomposition method \cite{Hoss06}, Variational iteration method \cite{Gha16}, Multi quadric method \cite{VanaA11}, Homotopy perturbation method \cite{Dehg10} and  Iterative schemes \cite{Nedia07}. But these methods have their limitations when dealing with nonlinear DAEs and sometimes complexity involved in calculations make them unsuitable for solving nonlinear DAEs.\smallskip

\noindent The regular form of DAEs is
\begin{equation} \label{eq1}
	G\left(w(v),w^{\prime}(v),v\right)= 0,  \hspace{1mm} G \in {C^1}({R^{2m + 1}},{R^m}), \hspace{1mm} v \in \left[0,V\right]
\end{equation}

\noindent where Jacobian $ \left[\frac{\partial G}{\partial w^{\prime}}\right] $  is singular on $R^{2m + 1}$ .\smallskip

\noindent Many DAEs arising in physical applications are in semi explicit form and while some other are in further restricted Hessenberg form \cite{Brenan}. \smallskip

\noindent The index-1 semi explicit DAEs are given by 
\begin{eqnarray} \label{eq3}
	\nonumber w^{\prime}(v) &=& G\left(w(v),u(v),v \right), G \in C({R^{m + k + 1}},{R^m}) \\
	0 &=& F\left(w(v),u(v),v\right),F \in C^{1}({R^{m + k + 1}},{R^k}),\hspace{1mm}v \in \left[0,V\right]
\end{eqnarray}
where $ \frac{\partial F}{\partial u }$ is nonsingular. \\

\noindent The index-2 Hessenberg DAEs are given by 
\begin{eqnarray} \label{eq5}
	\nonumber w^{\prime}(v) &=& G\left(w(v),u(v),v \right), G \in C^{1}({R^{m + k + 1}},{R^m})\\
	0 &=& F\left( w(v),v \right), F \in {C^2}({R^{m + 1}},{R^k}), \hspace{1mm} v \in \left[0,V\right]
\end{eqnarray}

\noindent where $\left(\frac{\partial F}{\partial u}\right)\left(\frac{\partial G}{\partial w} \right)$  is nonsingular \cite{Dehg10,Wang01}.\smallskip \smallskip

\noindent The index-3 Hessenberg DAEs are given by 
\begin{eqnarray} \label{eq5.1}
	\nonumber u^{\prime}(v) &=& G\left(w(v),u(v),s,v \right), G \in C^{1}({R^{m + k + 1}},{R^m}) \\
	\nonumber w^{\prime}(v) &=& F\left( w(v),u(v),v \right), F \in C^{2}({R^{m + k + 1}},{R^m}) \\
	0 &=& H\left( w(v),v\right), H \in {C^3}({R^{m + 1}},{R^l}),\hspace{1mm} v \in (0,V)
\end{eqnarray}

\noindent where $\left(\frac{\partial H}{\partial w}\right)\left(\frac{\partial F}{\partial u} \right)\left(\frac{\partial G}{\partial s}\right)$  is nonsingular \cite{Dehg10,Wang01}.
\smallskip \smallskip

The motivations for present work are the research work of authors \cite{Fatm04,Ben15,Ben16,Hoss06,VanaA11} who studied DAEs using various semi analytical methods, but these methods are not suitable to deal highly non-linear DAEs. Therefore some techniques are needed to overcome limitations of existing techniques and can directly solve nonlinear DAEs in  a well defined and reliable algorithm. 
\smallskip

\noindent We propose a simple approach involving the differential transformation in this paper. The differential transformation has been introduced by G. Pukhov as the ``Taylor transform" in 1976 and applied to the study of electrical circuits \cite{Puk82}. The differential transformation is closely related to Taylor expansion of real analytic functions. It has applications in solving different types of problems for all classes of differential equations (ordinary, partial, delayed, fractional, fuzzy etc.). The recent developments and applications of DTM are discussed in \cite{Aroz06,Bia10,Gok12,JafA11,JawH18,KaraB09,KangA09,Gir06,Gir19,Ravi08,Reb17,Re19,Reb20,Yang20} and references therein. \smallskip

\noindent In the present paper, the differential transformation is used to solve nonlinear differenial algebraic equations. The nonlinearity in the problems is addressed by using the partial ordinary Bell polynomials in the Fa\`{a} di Bruno's formula. The results obtained by this technique are compared to other methods. Error analysis of the method is studied for convergence criterion. To show the efficiency of the method some examples of this class are considered. These examples not only validate the accuracy of the method but also gives results which are more convergent to the exact solution. However, to the best of our knowledge, no researcher has applied the DTM using Bell polynomials on the practical problems discussed in the Section \ref{applications}.\smallskip

\noindent The paper is organized as follows. In Section \ref{preliminaries}, we introduce the main idea and basic formulae of the differential transformation and provide necessary results for the nonlinearities involving partial ordinary Bell polynomials. In Section \ref{estimation} we introduce the error estimate result. Numerical results and discussion are presented in Section \ref{applications}. A conclusion is given in section \ref{conclusion}. 

\section{Preliminaries}\label{preliminaries}
In this section we discuss the main idea and basic formulae of the differential transformation as well as notations and results related to transformation of general nonlinear terms.
\subsection{Idea of  differential transform}
Let $w(v)$ be analytical function in domain $D$ and $v=v_{0}$ be any arbitrary point in $D$. Then, $w(v)$  can be expaneded in series form about the point  $v=v_{0}$. 
The differential transform of the kth derivative of the function $w(v)$ is defined as
\begin{equation} \label{3.1} 
	W(k)[v_{0}]=\frac{1}{k!} \left[\frac{d^{k} w(v)}{dv^{k} } \right]_{v=v_{0}}.
\end{equation}

\noindent The inverse differential transformation is given by 
\begin{equation} \label{4.1} 
	w(v)=\sum _{k=0}^{\infty }W(k)[v_{0}](v-v_{0})^{k} .  
\end{equation} 

\noindent Using equation \eqref{3.1}-\eqref{4.1}

\begin{equation} \label{4.2} 
	w(v)=\sum _{k=0}^{\infty }\frac{1}{k!} \left[\frac{d^{k} w(v)}{dv^{k} } \right]_{v=v_{0}}(v-v_{0})^{k} . 
\end{equation} 

\noindent In real applications, the function $w(v)$ is expressed by finite sum

\begin{equation} \label{eq8} 
	w(v)=\sum _{k=0}^{N}W(k)[v_{0}](v-v_{0})^{k} .  
\end{equation} 

\noindent The results which are used in this paper are listed in table \eqref{tab} without proofs.

\begin{table}[htbp]
	\centering	
	\caption{Formulae of the differential transform method}
	\begin{tabular}{lll} \hline\noalign{\smallskip}
		& Original function & Transformed function \\ \noalign{\smallskip}\hline\noalign{\smallskip}
		1 & $\frac{d^{n} w(v)}{dv^{n} } $  & $(k+1)(k+2)(k+3){\dots}(k+n)W(k+n)$ \\ 
		2 & $w(v)=v^{n} $ &  $\delta (k-n)$, where $\delta \left(k - n\right) = \left\{ {\begin{array}{*{30}{c}}
				{1,{\rm{ }}k = n}\\
				{0,{\rm{ }}k \ne n}
		\end{array}} \right.$ \\ \vspace{1mm}
		3 &  $e^{\alpha v} $ &  $\frac{\alpha ^{k} }{k!} $	\\  \vspace{1mm}
		4 & $w_{1} (v)w_{2} (v)$       &   $\sum _{i=0}^{k} W_{1} (i)W_{2} (k-i) $\\ \noalign{\smallskip}\hline
	\end{tabular}
	\label{tab}
	
\end{table}

\subsection{Fa\`{a} di Bruno's formula and Bell polynomials}
One of the principal disadvantages of most papers based on applications of differential transformations is the differential transformation is not applied directly to nonlinear terms like $w^{n}, n\in {\mathbb N}$ or $e^{w}$. Authors \cite{Reb19} used Adomian polynomials to compute the differential transform of nonlinear terms.  However, the differential transformation of nonlinear terms can be determined without calculating and evaluating symbolic derivatives by applying Fa\`{a} di Bruno's formula to nonlinear terms.

\noindent Here we present some necessary notations and results obtained in \cite{Reb18}. The proofs are not included since they can be found in the cited paper.

\begin{definition}\cite{Comt74}
	The partial ordinary Bell polynomials are the polynomials $\hat{B}_{k,l}\left(\hat{x}_1,\dots,\hat{x}_{k-l+1} \right)$ in an infinite number of variables $\hat{x}_1, \hat{x}_2,\dots $ defined by the series expansion
	
	\begin{equation}\label{2.1}
		\sum\limits_{k \ge l} {{\hat{B}_{k,l}}\left( {\hat{x}_{1},\dots,\hat{x}_{k - l + 1}} \right)}t^{k}  = {\left( {\sum\limits_{m \ge 1} {{\hat{x}_m}} t^{m} } \right)^l},l = 0,1,2,\dots
	\end{equation}
\end{definition}

\begin{lemma}\cite{Reb18}\label{lemma1}
	The partial ordinary Bell polynomials $\hat{B}_{k,l}\left(\hat{x}_1,\dots,\hat{x}_{k-l+1} \right),l = 0,1,2,\dots,k\ge l$ satisfy the recurrence relation
	\begin{equation} \label{2.2}
		\hat{B}_{k,l}\left(\hat{x}_1,\dots,\hat{x}_{k-l+1} \right)= \sum\limits_{i = 1}^{k-l+1} \frac{i.l}{k}\hat{x}_{i}\hat{B}_{k-i,l-1}\left(\hat{x}_1,\dots,\hat{x}_{k-i-l+2} \right)
	\end{equation}
	where $\hat{B}_{0,0}=1$ and $\hat{B}_{k,0}=0$ for $k \ge 1$.
\end{lemma}

\begin{theorem}\cite{Reb18} \label{th2}
	Let $g$ and $f$ be real functions analytic near $t_{0}$ and $g(t_{0})$ respectively, and let $h$ be the composition $h\left(t \right)=\left({fog}\right)\left(t\right)=f\left( {g\left(t\right)}\right)$. Denote $D\left\{ {g\left( t \right)} \right\}\left[ {{t_0}} \right] = \left\{ {G\left( k \right)} \right\}_{k = 0}^\infty$, $D\left\{ {f\left( t \right)} \right\}\left[ {g\left( {{t_0}} \right)} \right] = \left\{ {F\left( k \right)} \right\}_{k = 0}^\infty $ and $D\left\{ {\left( {f \circ g} \right)\left( t \right)} \right\}\left[ {{t_0}} \right] = \left\{ {H\left( k \right)} \right\}_{k = 0}^\infty$   the differential transformations of functions $g$, $f$ and $h$ at $t_{0}$, $g\left({t_0} \right)$  and $t_{0}$   respectively. Then the numbers $H(k)$ in the sequence $\left\{ {H\left( k \right)} \right\}_{k = 0}^\infty$ satisfy the relations $ H(0)=F(0)$ and
	
	\begin{equation}
		H\left( k \right) = \sum\limits_{l = 1}^k {F\left( l \right).{{\hat B}_{k,l}}\left( {G\left( 1 \right),\dots,G\left( {k - l + 1} \right)} \right)}  \ \text{ for } k \geq 1.
	\end{equation} 
\end{theorem}

\subsection{Implementation of method}
Consider higher-index Hessenberg DAEs  as
\begin{eqnarray} \label{eq12.1}
	\nonumber w^{(m)}\left(v\right)&=&f\left(w(v), u\left(v\right)\right),  \\
	0&=&g\left(w(v)\right),
\end{eqnarray}
with initial conditions
\begin{equation} \label{eq12.2}
	w^{(i)}\left(0\right)=\eta_{i}, i=0,1,\dots,m-1,
\end{equation}
where $w^{(m)}$ is the $m^{th}$ derivatives of $w$ and $\eta_{i}$ are given constants.

\begin{figure}[htb]
	\centering
	{%
		\resizebox*{8cm}{!}{\includegraphics{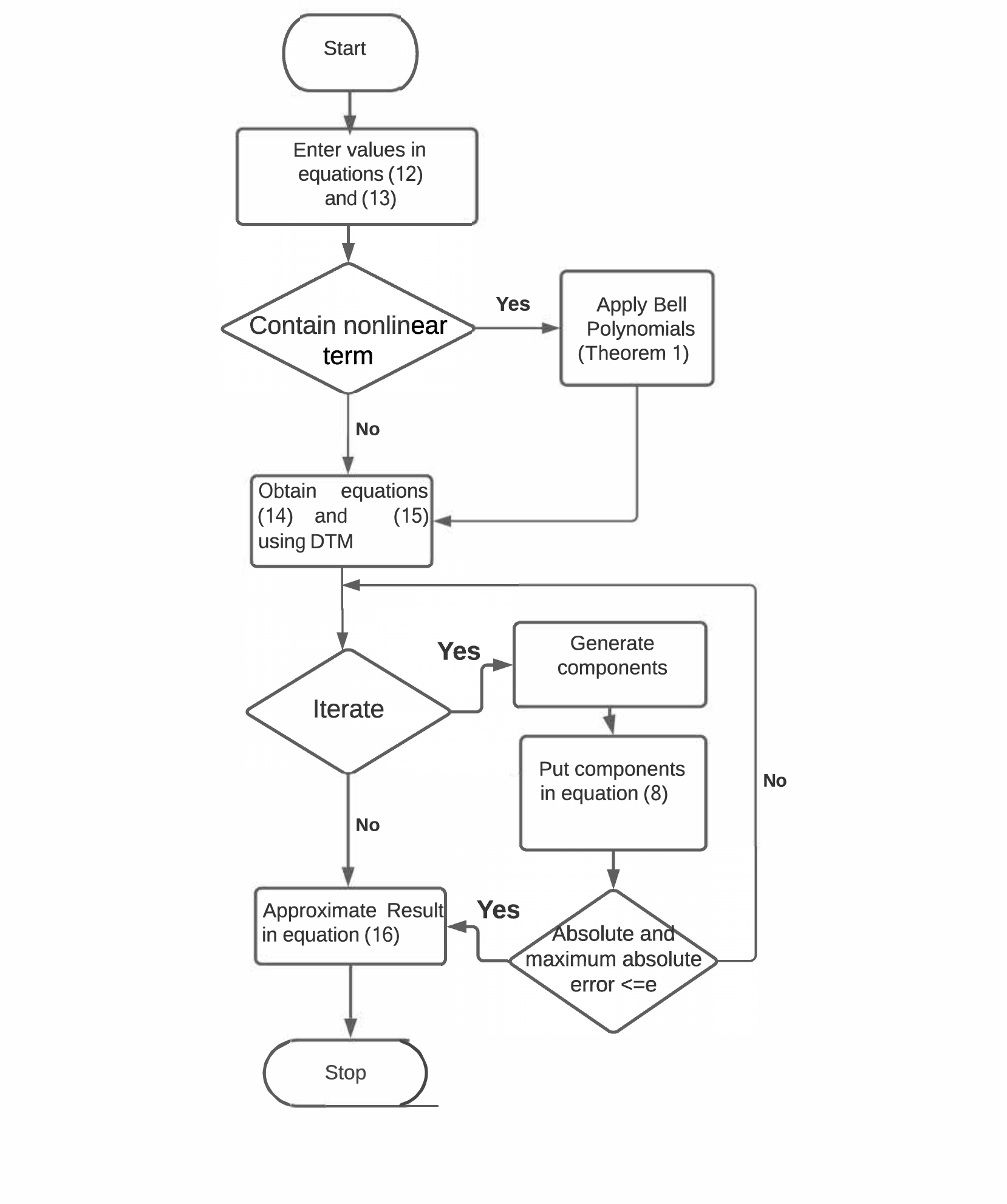}}}\hspace{5pt}
	\caption{Flow chart of present method.}
	\label{fig1}
\end{figure}

\noindent To solve equations \eqref{eq12.1} and \eqref{eq12.2}, apply differential transform we get the algebraic system
\begin{eqnarray} \label{eq13.1}
	\nonumber \left(k+1\right)\left(k+2\right)\dots \left(k+m\right)W\left(k+m\right)&=&F\left(W(k), U\left(k\right)\right),  \\
	0&=&G\left(W(k)\right),
\end{eqnarray}
and

\begin{equation} \label{eq13.2}
	W\left(k\right)=\eta_{i}, k=0,1,\dots,m-1,
\end{equation}
where $W$ and $U$  are the differential transform of $w$ and $u$ respectively and differential trasform of nonlinear term is obtained by Theorem \ref{th2}. The solution steps of equations \eqref{eq12.1}-\eqref{eq13.2} are explained in figure \ref{fig1}.

\noindent Now, the series solution is given by

\begin{equation} \label{14.1} 
	w(v)=\sum _{k=0}^{\infty }W(k)v^{k} .  
\end{equation} 

\begin{remark}
	Every step of the present method is illustrated through flowchart diagram shown in figure \ref{fig1} and implemented in example \ref{exam1} in section 4.
\end{remark}

\section{Error estimation}\label{estimation}

For comparison, absolute error and maximum absolute error are computed and defined as 
\begin{align*}
	E_N(v)&:=|w\left(v\right)-w_{N}\left(v\right)|, \\
	E_{N,\infty}&:=\mathop {\max}\limits_{0 \le v \le 1} E_N(v), 
\end{align*}
where $w\left(v\right)$  is the exact solution and $w_{N}\left(v\right)$  is the truncated series solution with degree $N$. Furthermore, the relative error between exact and approximate solution is defined by
$$
R_N(v):=\frac{E_N(v)}{|w(v)|}.
$$

Further following notations have been used in presented Tables:
\begin{eqnarray*}
	w_{i,N}(v)&:=&\text{Approximate solution obtained by present technique},\\
	w_{i,N}(v)\left[5\right]&:=&\text{Approximate solution obtained by Adomian method},\\
	w_{i,N}(v)\left[34\right]&:=&\text{Approximate solution obtained by Multi Quadric method},\\
	w_{i,N}(v)\left[21\right]&:=&\text{Approximate solution obtained by Lie Group method},\\
	R_{i,N}(v)&:=&\text{Relative error between exact and present solution.}
\end{eqnarray*}

\section{Numerical results and discussion}\label{applications}

Four examples of nonlinear higher index Hessenberg DAEs are solved to demonstrate the effectiveness of the proposed method. Mathematica software version 11 is used to perform all numerical computations. For reader's benefit every step of the proposed technique is explained in detail in example 1. 

\subsection{Example 1}\label{exam1}

\noindent Consider the nonlinear differential algebraic equation  \cite{Ben15}
\begin{eqnarray} \label{eq12} 
	\nonumber \frac{dw_{1}}{dv}&=&2w_{3}, \\
	\nonumber  \frac{dw_{2}}{dv}&=&2w_{4},\\
	\nonumber \frac{dw_{3}}{dv}&=&-2w_{3}+e^{w_{2}}+w+\phi_{1}, \\
	\nonumber \frac{dw_{4}}{dv}&=&2w_{4}+e^{w_{1}}+w+\phi_{2}, \\
	0&=&w_{1}+w_{2}-\phi_{3},\hspace{1mm} 0\le v<1,	
\end{eqnarray}
where
\begin{eqnarray}
	\phi_{1}\left(v\right)=-\frac{2v^{4}+2v^{3}+1}{2\left(1+v\right)^{2}},\hspace{1mm}
	\phi_{2}\left(v\right)=\frac{-2v^{4}+2v^{3}-1}{2\left(1-v\right)^{2}},\hspace{1mm}
	\phi_{3}\left(v\right)=ln\left(1-v^{2}\right),	\nonumber
\end{eqnarray}

\noindent initial conditions
\begin{equation} \label{eq13}
	w_{1} \left(0\right)=w_{2} \left(0\right)=0,\hspace{1mm} w_{3} \left(0\right)=\frac{1}{2},\hspace{1mm}w_{4} \left(0\right)=-\frac{1}{2}.
\end{equation}
The exact solution is given by 
\begin{eqnarray} \label{eq14} 
	\nonumber w_{1}\left(v\right)=ln\left(1+v\right),\hspace{1mm} w_{2}\left(v\right)=ln\left(1-v\right),\hspace{1mm}\\	w_{3}\left(v\right)=\frac{1}{2\left(1+v\right)},\hspace{1mm}w_{4}\left(v\right)=-\frac{1}{2\left(1-v\right)},\hspace{1mm}w\left(v\right)=v^{2}.
\end{eqnarray} 

Denoting 
\begin{eqnarray}
	\nonumber h_{1}\left( v \right) &=& f_{1}\left( g_{1}\left( v \right) \right), \text{where}\hspace{1mm} g_{1}\left( v \right) = w_{2}\left( v\right) \text{and}\hspace{1mm} f_{1}\left( x \right) = e^{x},	\\ \nonumber
	h_{2}\left( v \right) &=& f_{2}\left( g_{2}\left( v \right) \right), \text{where}\hspace{1mm} g_{2}\left( v \right) = w_{1}\left(v \right) \text{and}\hspace{1mm} f_{2}\left( x \right) = e^{x}.
\end{eqnarray}
Differential transformation of $f_{i}\left(x \right)$ is represented by $F_{i}\left( x \right)$  where $i=1,2$. Then, we get
\begin{equation}\label{eq15}
	F_{1}\left( k \right) = F_{2}\left( k \right) = \frac{1}{k!}.
\end{equation} 
The differential transform of $h_{1}\left(v\right)$ and $h_{2}\left(v\right)$ are represented using theorem \ref{th2}
\begin{eqnarray} \label{eq16}
	\nonumber H_{1}\left( 0 \right) &=& 1,\hspace{1mm} H_{1}\left(k \right) = \sum\limits_{l = 1}^k {{F_1}\left( l \right){{\hat B}_{k,l}}} \left( {W_{2}\left( 1 \right),\dots,W_{2}\left(k - l + 1 \right)} \right), \\
	H_{2}\left( 0 \right) &=& 1, \hspace{1mm} H_{2}\left( k \right) = \sum\limits_{l = 1}^k {F_2\left( l \right){{\hat B}_{k,l}}} \left( {W_{1}\left( 1 \right),\dots,W_{1}\left( k - l + 1 \right)} \right).	
\end{eqnarray}

Applying differential transform to equations \eqref{eq12}-\eqref{eq13}, we obtain the following recurrence relation
\begin{eqnarray} \label{eq17}
	\nonumber W_{1}\left(k + 1\right) &=& \frac{2}{\left(k + 1 \right)}W_{3}\left(k\right),\\ 
	\nonumber W_{2}\left(k + 1\right) &=& \frac{2}{\left(k + 1 \right)}W_{4}\left(k\right),\\ 
	\nonumber W_{3}\left(k + 1\right) &=& \frac{1}{\left(k + 1 \right)}\left(-2W_{3}\left(k\right)+H_{1}\left(k\right)+W\left(k\right)+\Phi_{1}\left(k\right)\right),\\ 
	\nonumber W_{4}\left(k + 1\right) &=& \frac{1}{\left(k + 1 \right)}\left(2W_{4}\left(k\right)+H_{2}\left(k\right)+W\left(k\right)+\Phi_{2}\left(k\right)\right),\\ 
	\nonumber 0&=&W_{1}\left(k\right)+W_{2}\left(k\right)-\Phi_{3}\left(k\right),\\
	W_{1}\left(0\right)&=&W_{2}\left(0\right)=0,\hspace{1mm} W_{3}\left(0\right)=\frac{1}{2},\hspace{1mm} W_{4}\left(0\right)=-\frac{1}{2},
\end{eqnarray}

where $\Phi_{1}$, $\Phi_{2}$ and $\Phi_{3}$ are differential transform of $\phi_{1}$, $\phi_{2}$ and $\phi_{3}$ respectively.

Using equations \eqref{eq15}-\eqref{eq17} we obtain the following components 
\begin{align} \label{eq23}
	\nonumber	k=0:\quad  W_{1}\left( 1 \right) &= 2W_{3}\left( 0 \right)=1, \\
	\nonumber W_{2}\left( 1 \right) &= 2W_{4}\left( 0 \right)=-1,\\ \nonumber
	W_{3}\left( 1 \right) &= -2W_{3}\left( 0 \right)+H_{1}\left(0\right)+W\left(0\right)+\Phi_{1}\left(0\right)=W\left(0\right)-\frac{1}{2},\\ \nonumber
	W_{4}\left( 1 \right) &= 2W_{4}\left( 0 \right)+H_{2}\left(0\right)+W\left(0\right)+\Phi_{2}\left(0\right)=W\left(0\right)-\frac{1}{2},\\ \nonumber
	\nonumber 0&=W_{1}\left(0\right)+W_{2}\left(0\right)-\Phi_{3}\left(0\right), \\
	\nonumber k=1:\quad  H_{1}\left( 1 \right) &=F_{1}\left( 1 \right){\hat B_{1,1}}\left( {W_{2}\left( 1 \right)} \right) = F_{1}\left( 1 \right)W_{2}\left( 1 \right) = -1, \\ \nonumber
	H_{2}\left( 1 \right) &=F_{2}\left( 1 \right){\hat B_{1,1}}\left( {W_{1}\left( 1 \right)} \right) = F_{2}\left( 1 \right)W_{1}\left( 1 \right) = 1, \\ \nonumber
	W_{1}\left(2 \right) &= W_{3}\left( 1\right), \\ \nonumber
	W_{2}\left( 2\right) &= W_{4}\left(1\right),\\ \nonumber
	W_{3}\left( 2\right) &=\frac{1}{2}\left( -2W_{3}\left(1 \right)+H_{1}\left(1\right)+W\left(1\right)+\Phi_{1}\left(1\right)\right)=\frac{1}{2}\left(-2W_{3}\left(1\right)+W\left(1\right)\right),\\ \nonumber
	W_{4}\left( 2\right) &=\frac{1}{2}\left( 2W_{4}\left(1 \right)+H_{2}\left(1\right)+W\left(1\right)+\Phi_{2}\left(1\right)\right)=\frac{1}{2}\left(2W_{4}\left(1\right)+W\left(1\right)\right),\\	  
	\nonumber 0&=W_{1}\left(1\right)+W_{2}\left(1\right)-\Phi_{3}\left(1\right),\\   
	\nonumber	k=2:\quad  H_{1}\left(2\right)&={\sum\limits_{l = 1}^2 {F_{1}\left( l \right){{\hat B}_{2,1}}} \left( W_{2}\left( 1 \right),W_{2}\left( 2 \right) \right)}= {F_{1}\left( 1 \right)W_{2}\left( 2 \right) + F_{1}\left( 2 \right){W_{2}^2}\left( 1 \right)}, \\ \nonumber
	H_{2}\left(2\right)&={\sum\limits_{l = 1}^2 {F_{2}\left( l \right){{\hat B}_{2,1}}} \left( W_{1}\left( 1 \right),W_{1}\left( 2 \right) \right)}= {F_{2}\left( 1 \right)W_{1}\left( 2 \right) + F_{2}\left( 2 \right){W_{1}^2}\left( 1 \right)}, \\ \nonumber
	W_{1}\left(3 \right) &=\frac{2}{3} W_{3}\left( 2\right), \\ \nonumber
	W_{2}\left( 3\right) &=\frac{2}{3} W_{4}\left(2\right),\\ \nonumber
	W_{3}\left( 3\right) &=\frac{1}{3}\left(-2W_{3}\left(2 \right)+H_{1}\left(2\right)+W\left(2\right)+\Phi_{1}\left(2\right)\right)=\frac{1}{3}\left(-2W_{3}\left(2\right)+W\left(2\right)-\frac{3}{2}\right),\\ \nonumber
	W_{4}\left( 3\right) &=\frac{1}{3}\left(2W_{4}\left(2 \right)+H_{2}\left(2\right)+W\left(2\right)+\Phi_{2}\left(2\right)\right)=\frac{1}{3}\left(2W_{4}\left(2\right)+W\left(2\right)-\frac{3}{2}\right),\\   
	0&=W_{1}\left(2\right)+W_{2}\left(2\right)-\Phi_{3}\left(2\right) \text{ and soon.}
\end{align}

Now, with the help of equation \eqref{eq8}, the series solution is given by
\begin{eqnarray} \label{eq46}
	\nonumber w_{1}\left( v \right) &=&  v - \frac{1}{2}{v^2}+\frac{1}{3}v^{3}-\frac{1}{4}v^{4}+ \dots, \\
	\nonumber w_{2}\left( v \right) &=& -v- \frac{1}{2}{v^2}-\frac{1}{3}v^{3}-\frac{1}{4}v^{4}- \dots,
\end{eqnarray}
\begin{eqnarray} 
	\nonumber w_{3}\left( v \right) &=&\frac{1}{2}-\frac{1}{2}v+ \frac{1}{2}{v^2}-\frac{1}{2}v^{3}+\frac{1}{2}v^{4}- \dots, \\
	\nonumber w_{4}\left( v \right) &=&-\frac{1}{2}-\frac{1}{2}v- \frac{1}{2}{v^2}-\frac{1}{2}v^{3}-\frac{1}{2}v^{4}- \dots, \\
	\nonumber w\left( v \right) &=& v^{2}, 
\end{eqnarray}
which converges to the exact solution given by equation \eqref{eq14}.

\begin{table}[htbp]
	\centering
	\caption{Comparison of numerical solution of $w_{1}(v)$ with exact solution and Benhammouda \cite{Ben15} solution  for example 1 ($N=20$)} 
	{\begin{tabular}{lllll} \hline\noalign{\smallskip}
			$v$ & $w_{1}\left(v\right)$ & $w_{1,N}(v)$ & $w_{1,N}(v)$ \cite{Ben15} & $R_{1,N}(v)$  \\ \noalign{\smallskip}\hline\noalign{\smallskip}
			0.1 & 0.0953101798 & 0.0953101798 & 0.0953101798   & 7.2E-16  \\ 
			0.2 & 0.1823215568  & 0.1823215568 & 0.1823215568 & 3.0E-16  \\ 
			0.3 & 0.2623642645 &0.2623642645 &0.2623642645   & 1.4E-12  \\ 
			0.4 & 0.3364722366 & 0.3364722365 & 0.3364722365 & 4.5E-10 \\ 
			0.5 & 0.4054651081 & 0.4054650927 & 0.4054650927 & 3.7E-08  \\ 
			0.6 & 0.4700036292 & 0.4700029649 & 0.4700029649& 1.4E-06  \\
			0.7 & 0.5306282511 & 0.5306123016 & 0.5306123016  & 3.0E-05  \\ 
			0.8 & 0.5877866649 & 0.5875375291&0.5875375291  & 4.2E-04  \\ 
			0.9 & 0.6418538862& 0.6390499221 & 0.6390499221  & 4.3E-03 \\ 
			1.0 & 0.6931471806 &0.6687714032 & 0.6687714032 & 3.5E-02  \\ \noalign{\smallskip}\hline
	\end{tabular}}
	\label{tab1.1}
\end{table}

\begin{table}[htbp]
	\centering
	\caption{Comparison of numerical solution of $w_{2}(v)$ with exact solution and Benhammouda \cite{Ben15} solution  for example 1 ($N=20$)}
	{\begin{tabular}{lllll} \hline\noalign{\smallskip}
			$v$ & $w_{2}\left(v\right)$ &$w_{2,N}\left(v\right)$  &  $w_{2,N}\left(v\right)$\cite{Ben15}  &  $R_{2,N}(v)$  \\ 	\noalign{\smallskip}\hline\noalign{\smallskip}
			0.1 & -0.1053605157 & -0.1053605157 & -0.1053605157   & 3.9E-16  \\ 
			0.2 & -0.2231435513 & -0.2231435513 & -0.2231435513  & 2.4E-16  \\ 
			0.3 & -0.3566749439 & -0.3566749439 & -0.3566749439  & 1.9E-12  \\ 
			0.4 & -0.5108256238 & -0.5108256234 & -0.5108256234  & 6.6E-10 \\ 
			0.5 & -0.6931471806 & -0.6931471371 & -0.6931471371 & 6.2E-08  \\ 
			0.6 & -0.9162907319 & -0.9162882787 & -0.9162882787 & 2.6E-06  \\
			0.7 & -1.2039728040 & -1.2038920520 & -1.2038920520 & 6.7E-05  \\ 
			0.8 & -1.6094379120 & -1.6075458670 & -1.6075458670   & 1.1E-03  \\ 
			0.9 & -2.3025850930 & -2.2633497340 & -2.2633497340 & 1.7E-02 \\  	\noalign{\smallskip}\hline
	\end{tabular}}
	\label{tab1.2}
\end{table}

\begin{table}[htbp]
	\centering
	\caption{Comparison of numerical solution of $w_{3}(v)$ with exact solution and Benhammouda \cite{Ben15} solution  for example 1(N=20)}
	{\begin{tabular}{lllll} \hline\noalign{\smallskip}
			$v$ & $w_{3}\left(v\right)$ &$w_{3,N}\left(v\right)$  &  $w_{3,N}\left(v\right)$\cite{Ben15}  &  $R_{3,N}(v)$  \\ 	\noalign{\smallskip}\hline\noalign{\smallskip}
			0.1 & 0.4545454545  & 0.4545454545 &0.4545454545  & 0  \\ 
			0.2 & 0.4166666667  & 0.4166666667 &0.4166666667   & 2.1E-15  \\ 
			0.3 & 0.3846153846  & 0.3846153846 &0.3846153846  & 1.0E-11  \\ 
			0.4 & 0.3571428571  & 0.3571428587 &0.3571428587  & 4.3E-09 \\ 
			0.5 & 0.3333333333  & 0.3333334923 &0.3333334923  & 4.7E-07  \\ 
			0.6 & 0.3125000000  & 0.3125068553 &0.3125068553  & 2.1E-05  \\
			0.7 & 0.2941176471  & 0.2942819253 &0.2942819253 & 5.5E-04  \\ 
			0.8 & 0.2777777778  & 0.2803398256 &0.2803398256 & 9.2E-03  \\ 
			0.9 & 0.2631578947  & 0.2919523656 &0.2919523656 & 1.0E-01 \\  \noalign{\smallskip}\hline
	\end{tabular}}
	\label{tab1.3}
\end{table}

\begin{table}[htbp]
	\centering
	\caption{Comparison of numerical solution of $w_{4}(v)$ with exact solution and Benhammouda \cite{Ben15} solution  for example 1(N=20)}
	{\begin{tabular}{lllll}	\hline\noalign{\smallskip}
			$v$ & $w_{4}\left(v\right)$ &$w_{4,N}\left(v\right)$  &  $w_{4,N}\left(v\right)$\cite{Ben15}  &  $R_{4,N}(v)$   \\ \noalign{\smallskip}\hline\noalign{\smallskip}
			0.1 & -0.5555561000 & -0.5555561000&-0.5555561000  & 0  \\ 
			0.2 & -0.6251000000 & -0.6251000000&-0.6251000000 & 2.1E-15  \\ 
			0.3 & -0.7142861000 & -0.7142861000&-0.7142861000  & 1.0E-11  \\ 
			0.4 & -0.8333331000 & -0.8333331000&-0.8333331000   & 4.3E-09 \\ 
			0.5 & -1.1000000000 & -1.1000000000&-1.1000000000 & 4.7E-07  \\ 
			0.6 & -1.2510000000 & -1.2499710000&-1.2499710000  & 2.1E-05  \\
			0.7 & -1.6666710000 & -1.6657410000&-1.6657410000 & 5.5E-04  \\ 
			0.8 & -2.5100000000 & -2.4769410000&-2.4769410000  & 9.2E-03  \\ 
			0.9 & -5.1000000000 & -4.4529110000&-4.4529110000 & 1.0E-01 \\  \noalign{\smallskip}\hline
	\end{tabular}}
	\label{tab1.4}
\end{table}

\begin{table}[htbp]
	\centering
	\caption{Maximum absolute error for $w_{1}$, $w_{2}$, $w_{3}$ and $w_{4}$ of example 1}
	{\begin{tabular}{lcccc} \hline\noalign{\smallskip}
			$N$ & $E_{1N,\infty}$& $E_{2N,\infty}$ & $E_{3N,\infty}$&$E_{4N,\infty}$ \\ \noalign{\smallskip}\hline\noalign{\smallskip}
			10 & 1.5E-02 & 1.8E-01 &  8.2E-02 &1.5E-00\\ 
			15 & 6.2E-03 & 8.27E-02 & 4.8E-02 &9.2E-01\\
			20 & 2.8E-03 & 3.9E-02 &  2.8E-02 &5.4E-01\\   \noalign{\smallskip}\hline
	\end{tabular}}
	\label{tab1.5}
\end{table}

\begin{figure}[H]
	\centering
	\subfloat[absolute error for $w_{1}(v)$]{%
		\resizebox*{6.5cm}{!}{\includegraphics{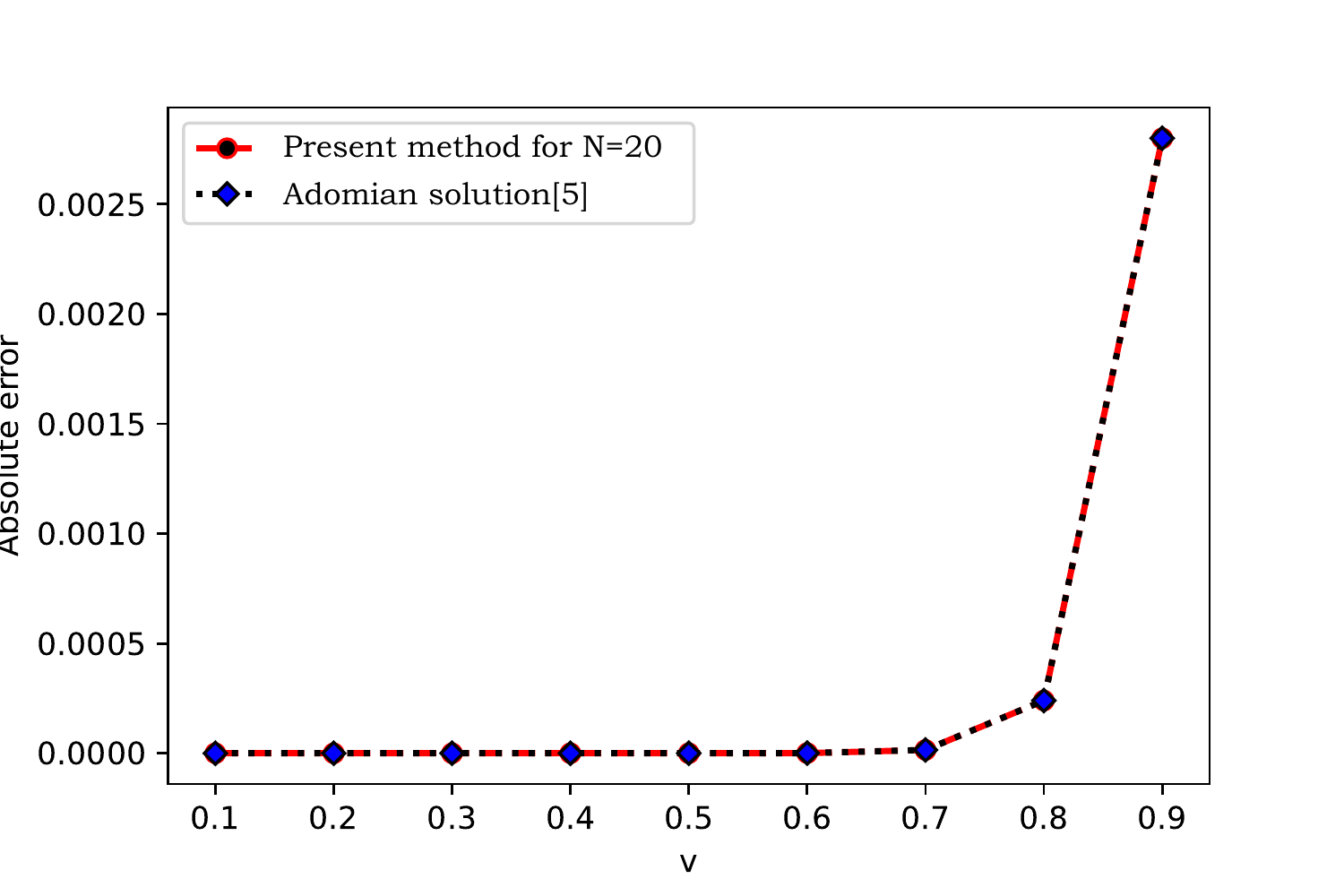}}}\hspace{2pt}
	\subfloat[absolute error for $w_{2}(v)$]{%
		\resizebox*{6.5cm}{!}{\includegraphics{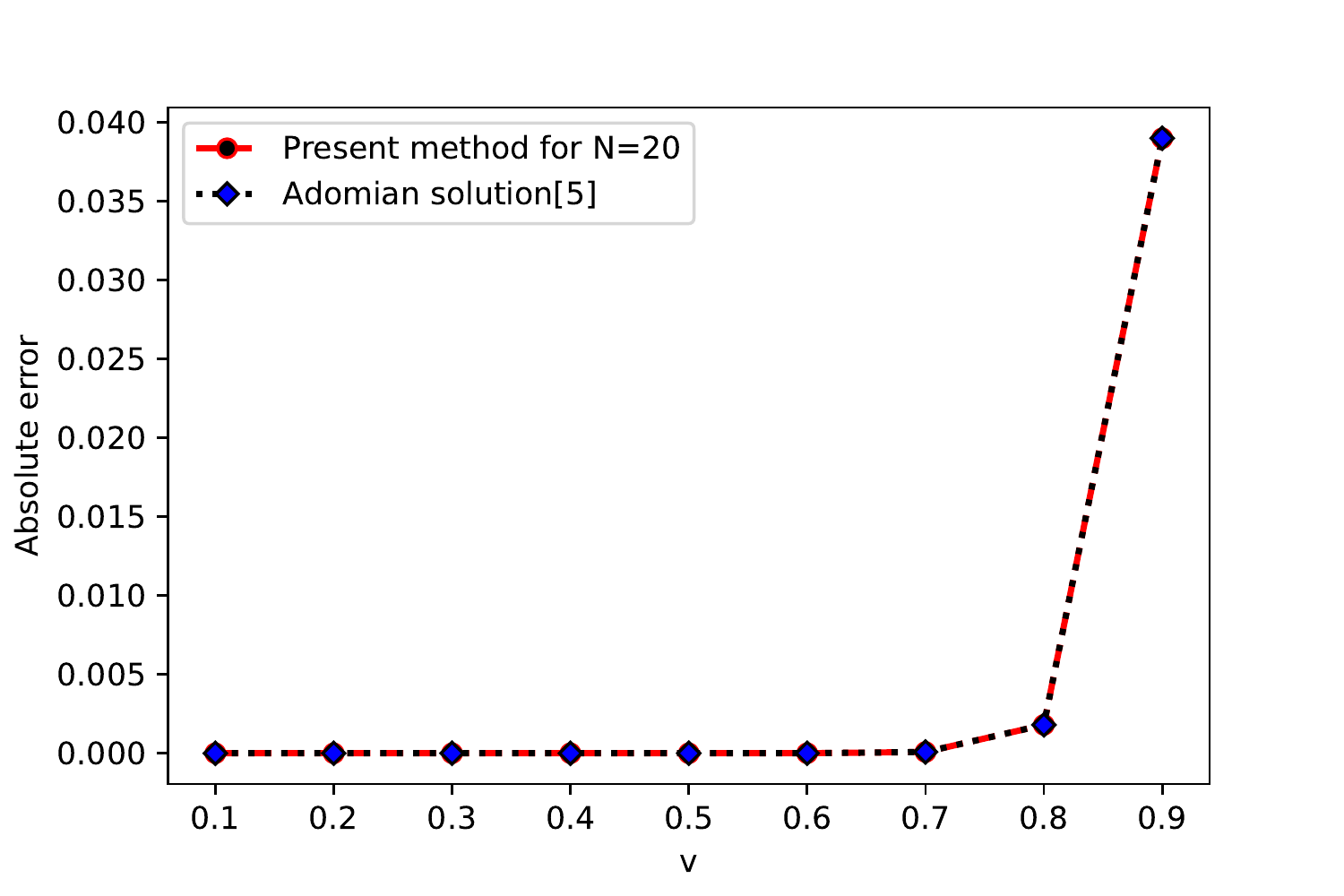}}}\hspace{2pt}
	\subfloat[absolute error for $w_{3}(v)$]{%
		\resizebox*{6.5cm}{!}{\includegraphics{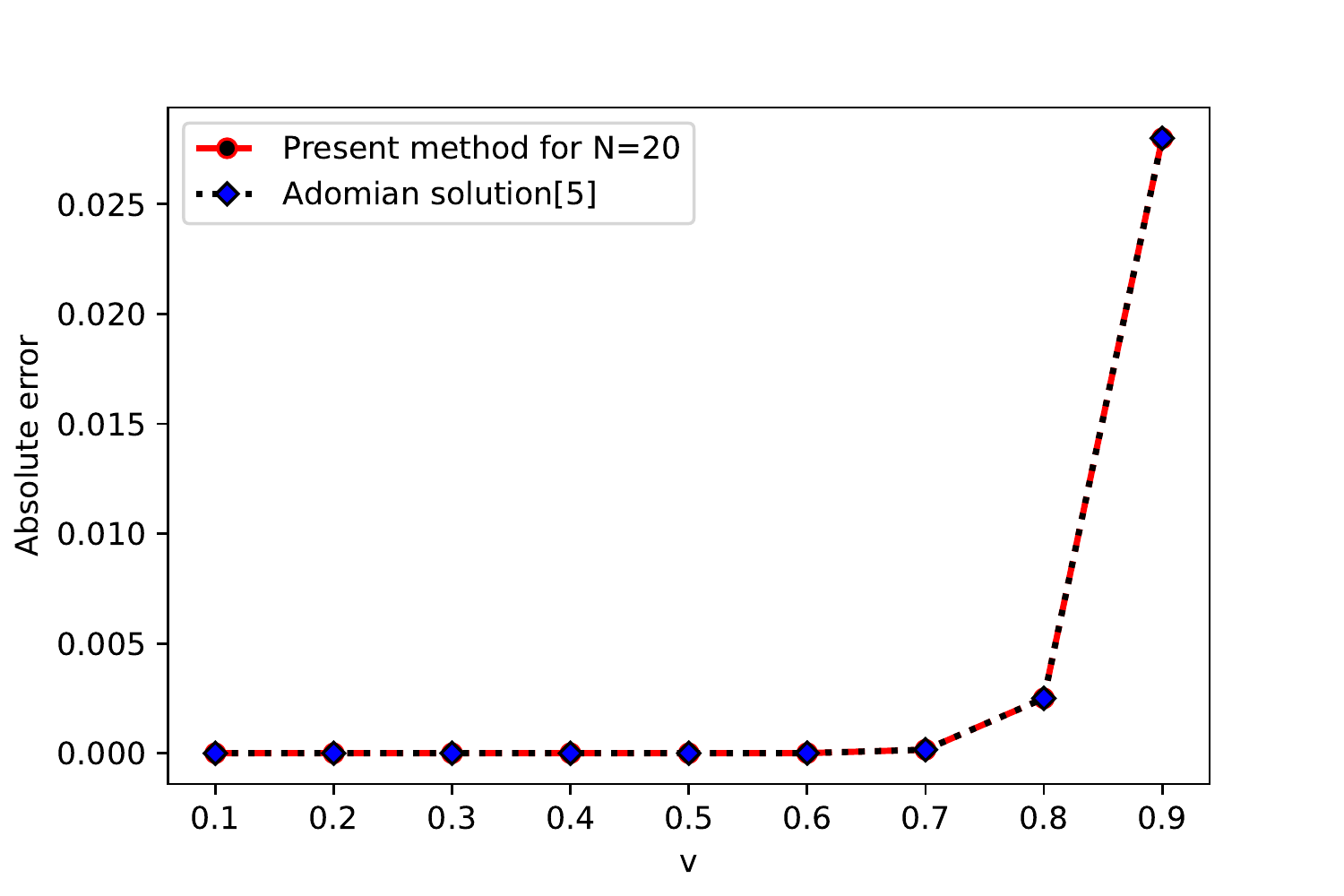}}}\hspace{2pt}
	\subfloat[absolute error for $w_{4}(v)$]{%
		\resizebox*{6.5cm}{!}{\includegraphics{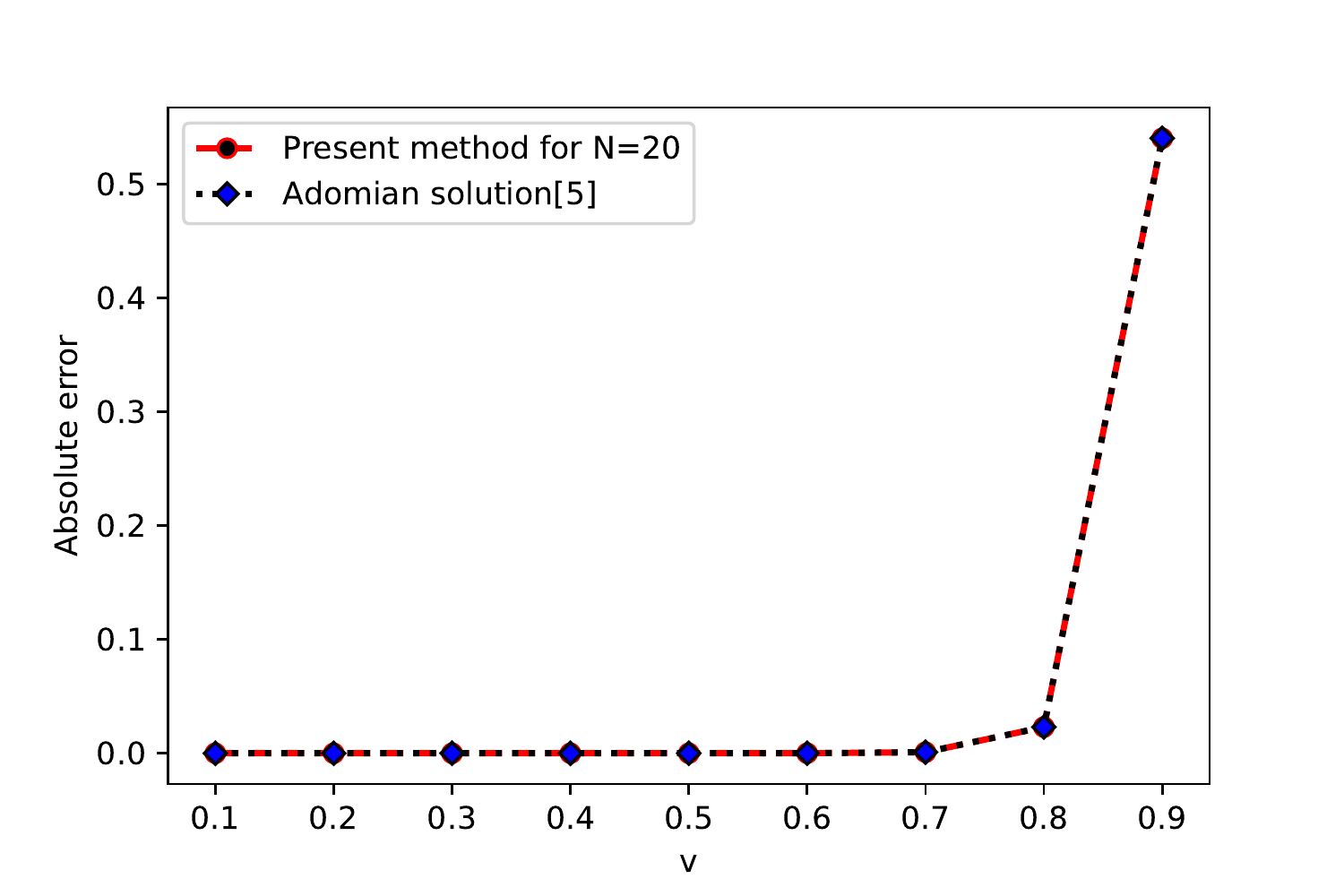}}}\hspace{2pt}
	\caption{The absolute error for $w_{1}$, $w_{2}$, $w_{3}$ and $w_{4}$ of example 1.} 
	\label{fig1.1}
\end{figure}

Tables [\eqref{tab1.1}, \eqref{tab1.2},\eqref{tab1.3},\eqref{tab1.4}] compare the proposed solution to the exact and numerical solutions presented in \cite{Ben15} for the variables $w_{1}$, $w_{2}$, $w_{3}$, and $w_{4}$. It is evident that the present solution conforms well with the exact solution and the solution given in \cite{Ben15}. In addition, the maximum absolute error is determined in table \eqref{tab1.5}, and it decreases as the number of series terms increases.

\subsection{Example 2}
\noindent Consider the nonlinear differential algebraic equation  \cite{Ben15}

\begin{eqnarray}
	\nonumber	\frac{d^{2}w_{1}}{dv^{2}}&=&2w_{2}-2w^{3}_{2}-w_{1}w,\\
	\nonumber	\frac{d^{2}w_{2}}{dv^{2}}&=&2w_{1}-2w^{3}_{1}-w_{2}w,\\
	0&=&w^{2}_{1}+w^{2}_{2}-1,\hspace{1mm} v \ge 0,
\end{eqnarray}

initial conditions
\begin{eqnarray}
	w_{1}\left(0\right)=1,\hspace{1mm} \frac{dw_{1}\left(0\right)}{dv}=0,\hspace{1mm} w_{2}\left(0\right)=0, \hspace{1mm}\frac{dw_{2}\left(0\right)}{dv}=1.
\end{eqnarray}

The exact solution is given by
\begin{eqnarray}
	w_{1}\left(v\right)=cosv,\hspace{1mm} w_{2}\left(v\right)=sinv,\hspace{1mm} w\left(v\right)=1+sin(2v).
\end{eqnarray}

\begin{table}[htbp]
	\centering
	\caption{Comparison of numerical solution of $w_{1}(v)$ with exact solution and Benhammouda \cite{Ben15} solution  for example 2 $\left(N=15\right)$}
	{\begin{tabular}{lllll} 	\hline\noalign{\smallskip}
			$v$ & $w_{1}\left(v\right)$ &$w_{1,N}\left(v\right)$  &  $w_{1,N}\left(v\right)$\cite{Ben15}  &  $R_{1,N}(v)$ \\ \noalign{\smallskip}\hline\noalign{\smallskip}
			0.1 &0.9950041000 &0.9950041000&0.9950041000 & 0  \\ 
			0.2 &0.9800671000 &0.9800671000&0.9800671000 & 0  \\ 
			0.3 &0.9553361000 &0.9553361000&0.9553361000  & 0  \\ 
			0.4 &0.9210611000 &0.9210611000&0.9210611000 & 1.2E-16 \\ 
			0.5 &0.8775831000 &0.8775831000&0.8775831000  & 1.2E-16  \\ 
			0.6 &0.8253361000 &0.8253361000&0.8253361000 & 1.2E-16  \\
			0.7 &0.7648421000 &0.7648421000&0.7648421000  & 2.9E-16  \\ 
			0.8 &0.6967071000 &0.6967071000&0.6967071000  & 2.0E-15  \\ 
			0.9 &0.6216110000 &0.6216110000&0.6216110000 & 1.4E-14\\  \noalign{\smallskip}\hline
	\end{tabular}}
	\label{tab2.1}
\end{table}

\begin{table}[htbp]
	\centering
	\caption{Comparison of numerical solution of $w_{2}(v)$ with exact solution and Benhammouda \cite{Ben15} solution  for example 2 $\left(N=15\right)$}
	{\begin{tabular}{lllll} 	\hline\noalign{\smallskip}
			$v$ & $w_{2}\left(v\right)$ &$w_{2,N}\left(v\right)$  &  $w_{2,N}\left(v\right)$\cite{Ben15}  &  $R_{2,N}(v)$  \\ \noalign{\smallskip}\hline\noalign{\smallskip}
			0.1 &0.0998334100 & 0.0998334100&0.0998334100 & 1.3E-16  \\ 
			0.2 &0.1986691000 & 0.1986691000&0.1986691000 & 1.3E-16  \\ 
			0.3 &0.2955210000 & 0.2955210000&0.2955210000  & 1.3E-16  \\ 
			0.4 &0.3894181000 & 0.3894181000&0.3894181000  & 1.3E-16 \\ 
			0.5 &0.4794261000 & 0.4794261000&0.4794261000   & 1.3E-16  \\ 
			0.6 &0.5646421000 & 0.5646421000&0.5646421000  & 1.3E-16  \\
			0.7 &0.6442181000 & 0.6442181000&0.6442181000 & 1.7E-16  \\ 
			0.8 & 0.717356100 & 0.7173561000&0.7173561000  & 1.5E-16  \\ 
			0.9 & 0.783327100 & 0.7833271000&0.7833271000 & 5.6E-16 \\  \noalign{\smallskip}\hline
	\end{tabular}}
	\label{tab2.2}
\end{table}

\begin{table}[htbp]
	\centering
	\caption{Comparison of numerical solution of $w(v)$ with exact solution and Benhammouda \cite{Ben15} solution  for example 2 $\left(N=15\right)$}
	{\begin{tabular}{lllll} \hline\noalign{\smallskip}
			$v$ & $w\left(v\right)$ &$w_{N}\left(v\right)$  &  $w_{N}\left(v\right)$\cite{Ben15}  &  $R_{N}(v)$  \\ 	\noalign{\smallskip}\hline\noalign{\smallskip}
			0.1 &1.1986710000 & 1.1986710000&1.1986710000 & 1.8E-16  \\ 
			0.2 &1.3894210000 & 1.3894210000&1.3894210000  & 1.8E-16  \\ 
			0.3 &1.5646410000 & 1.5646410000&1.5646410000  & 1.8E-16  \\ 
			0.4 &1.7173610000 & 1.7173610000&1.7173610000  & 1.2E-16 \\ 
			0.5 &1.8414710000 & 1.8414710000&1.8414710000   & 1.5E-15  \\ 
			0.6 &1.9320410000 & 1.9320410000&1.9320410000  & 3.2E-14  \\
			0.7 &1.9854510000 & 1.9854510000&1.9854510000 & 4.2E-13  \\ 
			0.8 &1.9995710000 & 1.9995710000&1.9995710000  & 4.1E-12  \\ 
			0.9 &1.9738510000 & 1.9738510000&1.9738510000 & 3.0E-11 \\  	\noalign{\smallskip}\hline
	\end{tabular}}
	\label{tab2.3}
\end{table}

\begin{figure}[htbp]
	\centering
	\subfloat[absolute error for $w_{1}(v)$]{%
		\resizebox*{6.8cm}{!}{\includegraphics{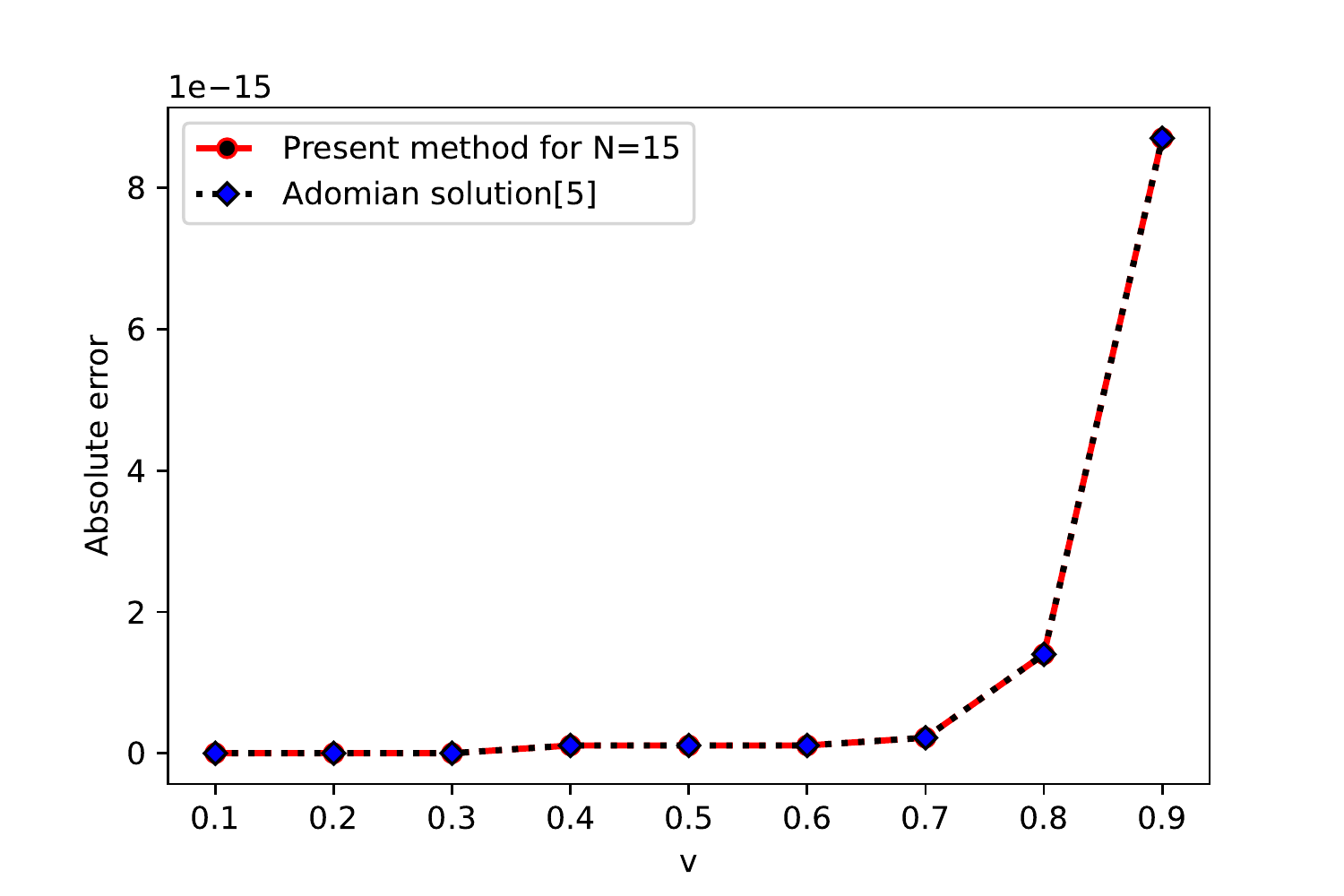}}}\hspace{3pt}
	\subfloat[absolute error for $w_{2}(v)$]{%
		\resizebox*{6.8cm}{!}{\includegraphics{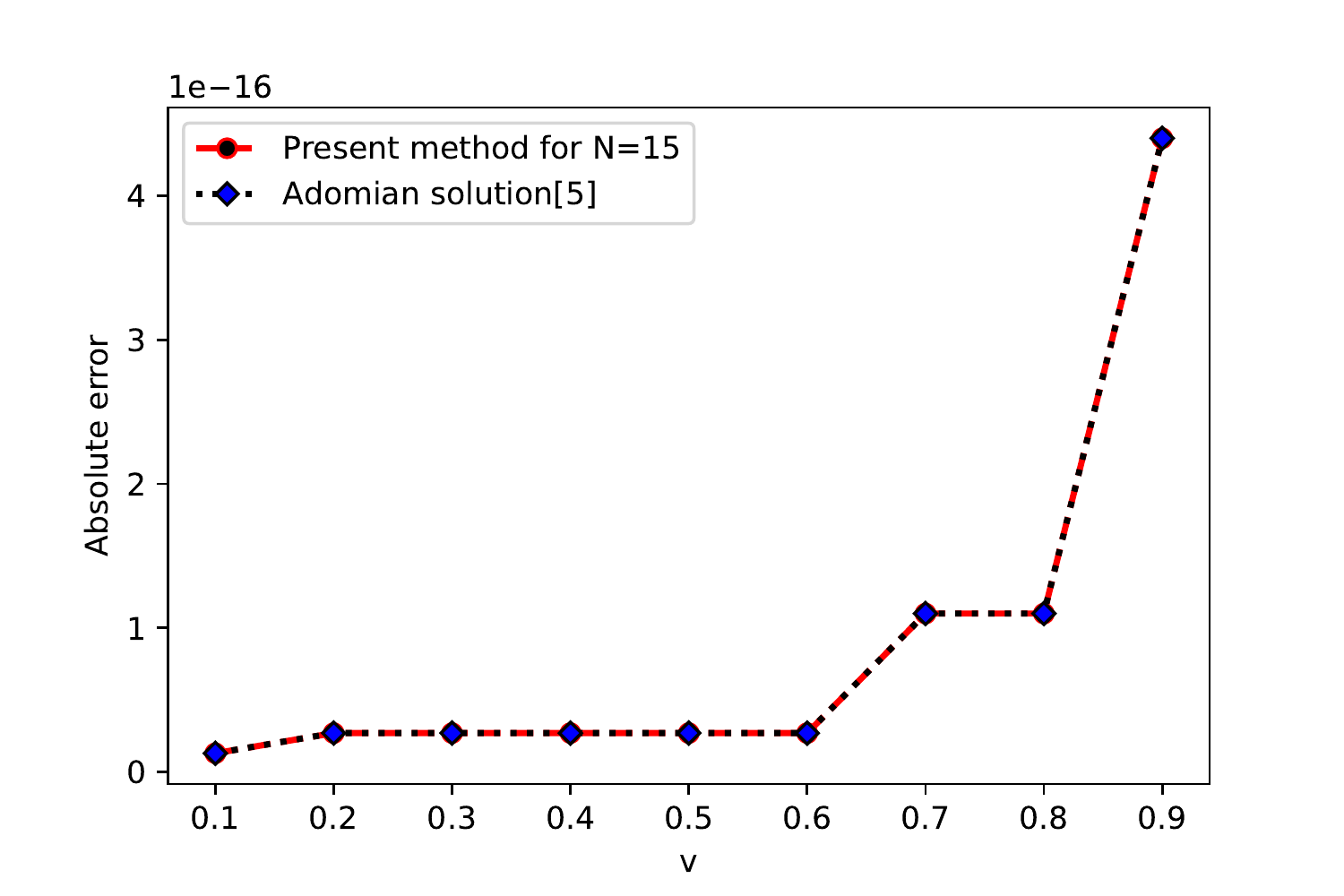}}}\hspace{3pt}
	\subfloat[absolute error for $w(v)$]{%
		\resizebox*{6.8cm}{!}{\includegraphics{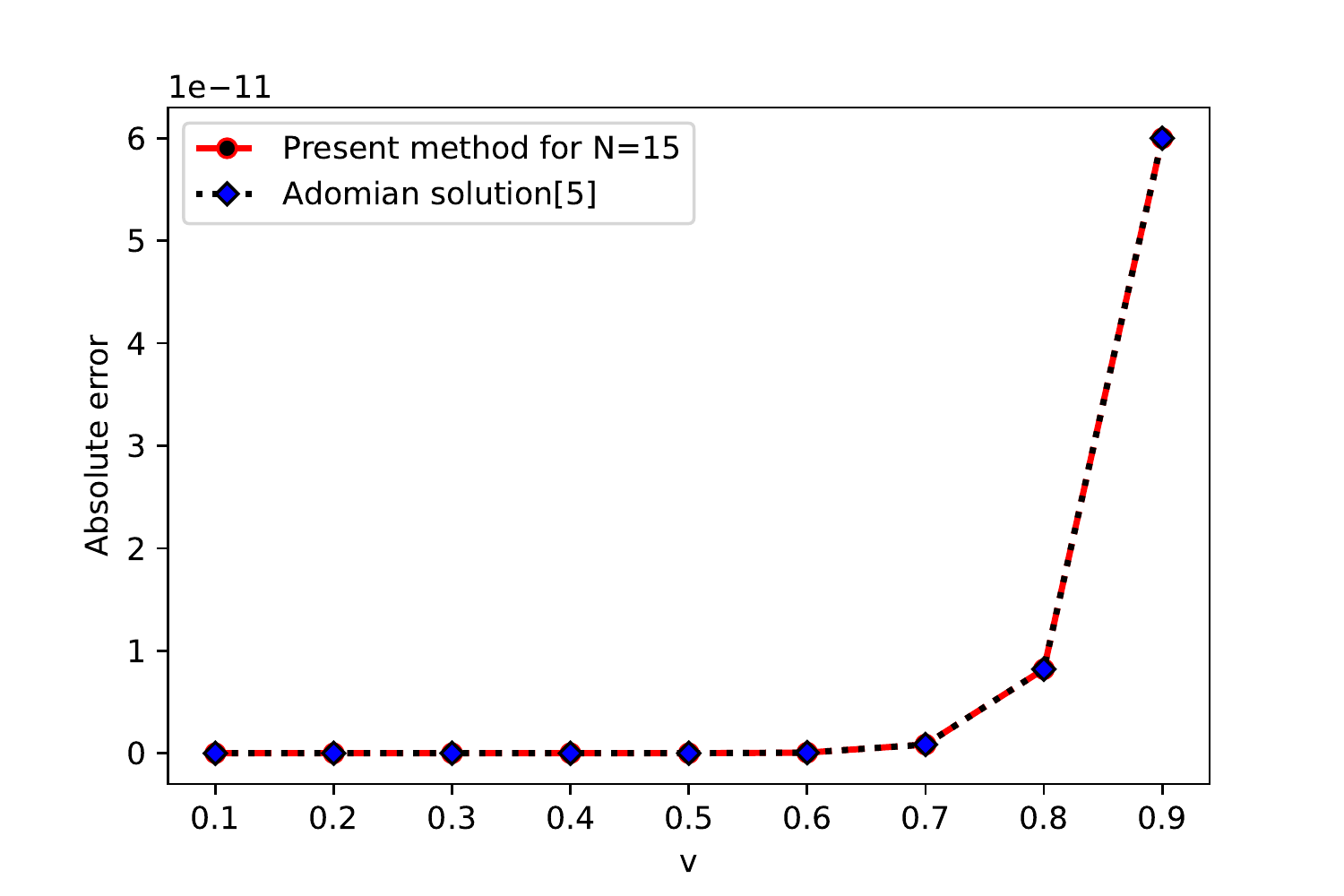}}}
	\caption{The absolute error for $w_{1}$, $w_{2}$ and $w$ of example 2.} 
	\label{fig2.1}
\end{figure}

\begin{table}[htbp]
	\centering
	\caption{Maximum absolute error for $w_{1}$, $w_{2}$ and $w$ of example 2}
	{\begin{tabular}{llll} \hline\noalign{\smallskip}
			$N$ &   $E_{1N,\infty}$ & $E_{2N,\infty}$&  $E_{N,\infty}$ \\ \noalign{\smallskip}\hline\noalign{\smallskip}
			5 & 7.2E-04 & 9.3E-05 &  1.1E-02 \\ 
			10 & 5.8E-10& 7.8E-09 &  1.5E-05 \\
			15 & 8.7E-15& 4.4E-16 &  6.0E-11 \\   	\noalign{\smallskip}\hline
	\end{tabular}}
	\label{tab2.4}
\end{table}

The present solution is compared to the exact solution and numerical solution discussed in \cite{Ben15}  for $w_{1}$, $w_{2}$, and $w$ in tables [\eqref{tab2.1}, \eqref{tab2.2},\eqref{tab2.3}] respectively. The present solution is clearly in good conformity with the exact solution and solution discussed in \cite{Ben15}. The maximum absolute error decreases as the number of terms in the series solutions increases, as shown in table \eqref{tab2.4}.

\subsection{Example 3}

\noindent Consider the nonlinear differential algebraic equation  \cite{VanaA11}
\begin{eqnarray} \label{1} 
	\nonumber \frac{dw_{1}}{dv}&=&v\lambda \frac{dw_{3}}{dv}+e^{v} -v \lambda \left(e^{v} +e^{-v}\right), \\
	\nonumber  \frac{dw_{2}}{dv}&=&\left(\lambda -5\right)\frac{dw_{3}}{dv}-e^{-v} -\left(\lambda -5\right)\left(e^{v} +e^{-v} \right), \\
	0&=&v^{2} w_{1} \left(v\right)+w_{2} \sin{v} -v^{2} e^{v} -\sin{v} e^{-v}, \hspace{1mm} 0 \le v \leq 1,
\end{eqnarray}

\noindent where $\lambda$ is arbitrary parameter and we take $\lambda=15$ and

\noindent initial conditions
\begin{equation} \label{4}
	w_{1} \left(0\right)=1,\hspace{1mm} w_{2} \left(0\right)=1,\hspace{1mm} w_{3}\left(0\right)=0.
\end{equation}
The exact solution is given by 
\begin{equation} \label{5} 
	w\left(v\right)=\left(\begin{array}{c} {e^{v} } \\ {e^{-v} } \\ {e^{v} -e^{-v} } \end{array}\right).
\end{equation} 

\begin{table}[htbp]
	\centering
	\caption{Comparison of numerical solution of $w_{1}(v)$ with exact solution and Vanani \cite{VanaA11} solution  for example 3 $\left(N=12\right)$}
	{\begin{tabular}{llllp{2.5cm}} 	\hline\noalign{\smallskip}
			$v$ & $w_{1}\left(v\right)$ &$w_{1,N}\left(v\right)$  &  $w_{1,N}\left(v\right)$\cite{VanaA11}  &  $R_{1,N}(v)$ \\ 	\noalign{\smallskip}\hline\noalign{\smallskip}
			0.0 & 1.000000000 & 1.000000000 &1.000000000  & 0  \\ 
			0.1 & 1.105170918 & 1.105170918 & 1.105170918  & 2.0E-16  \\ 
			0.2 & 1.221402758 & 1.221402758 &1.221402758  & 2.0E-16  \\ 
			0.3 & 1.349858807 & 1.349858807 & 1.349858807  & 1.6E-16  \\ 
			0.4 & 1.491824697 & 1.491824697 &1.491824697  & 7.4E-16 \\ 
			0.5 & 1.648721270 & 1.648721270 &1.648721270  & 1.2E-14  \\ 
			0.6 & 1.822118800 & 1.822118800 &1.822118800  & 1.2E-13  \\
			0.7 & 2.013752707 & 2.013752707 & 2.013752707 & 8.1E-13  \\ 
			0.8 & 2.225540928 & 2.225540928 &2.225540928  & 4.2E-12  \\ 
			0.9 & 2.459603111 & 2.459603111 & 2.459603111 & 1.7E-11 \\ 
			1.0 & 2.718281828 & 2.718281828 &2.718281828  & 6.3E-11  \\ 	\noalign{\smallskip}\hline
	\end{tabular}}
	\label{tab1}
\end{table}

\begin{table}[htbp]
	\centering
	\caption{Comparison of numerical solution of $w_{2}(v)$ with exact solution and Vanani \cite{VanaA11} solution  for example 3 $\left(N=12\right)$}
	{\begin{tabular}{lllll} \hline\noalign{\smallskip}
			$v$ & $w_{2}\left(v\right)$ &$w_{2,N}\left(v\right)$  &  $w_{2,N}\left(v\right)$\cite{VanaA11}  &  $R_{2,N}(v)$ \\ \noalign{\smallskip}\hline\noalign{\smallskip}
			0.0 & 1.000000000 & 1.000000000 &1.000000000   & 0  \\ 
			0.1 & 0.904837418 & 0.904837418 & 0.904837418  & 1.2E-16  \\ 
			0.2 & 0.818730753 & 0.818730753 & 0.818730753  & 1.3E-16  \\ 
			0.3 & 0.740818220 & 0.740818220 & 0.740818220  & 1.3E-16  \\ 
			0.4 & 0.670320046 & 0.670320046 & 0.670320046  & 1.4E-15 \\ 
			0.5 & 0.606530659 & 0.606530659 & 0.606530659   & 3.1E-14  \\ 
			0.6 & 0.548811636 & 0.548811636 & 0.548811636  & 3.6E-13  \\
			0.7 & 0.496585303 & 0.496585303 & 0.496585303  & 2.9E-12  \\ 
			0.8 & 0.449328964 & 0.449328964 & 0.449328964  & 1.8E-11  \\ 
			0.9 & 0.406569659 & 0.406569659 & 0.406569659  & 9.4E-11 \\ 
			1.0 & 0.367879441 & 0.367879441 & 0.367879441   & 4.0E-10  \\ \noalign{\smallskip}\hline
	\end{tabular}}
	\label{tab2}
\end{table}

\begin{table}[htbp]
	\centering
	\caption{Comparison of numerical solution of $w_{3}(v)$ with exact solution and Vanani \cite{VanaA11} solution  for example 3 $\left(N=12\right)$}
	{\begin{tabular}{lllll} \hline\noalign{\smallskip}
			$v$ & $w_{3}\left(v\right)$ &$w_{3,N}\left(v\right)$  &  $w_{3,N}\left(v\right)$\cite{VanaA11}  &  $R_{3,N}(v)$ \\ \noalign{\smallskip}\hline\noalign{\smallskip}
			0.0 & 0.000000000 & 0.000000000 &0.000000000  & 0  \\ 
			0.1 & 0.200333500 & 0.200333500 & 0.200333500   & 6.9E-16  \\ 
			0.2 & 0.402672005 & 0.402672005 & 0.402672005 & 1.3E-16  \\ 
			0.3 & 0.609040586 & 0.609040586&0.609040586  & 1.3E-16  \\ 
			0.4 & 0.821504651 & 0.821504651 & 0.821504651  & 2.5E-15 \\ 
			0.5 & 1.042190610 & 1.042190610 &1.042190610 & 3.7E-14  \\ 
			0.6 & 1.273307164 & 1.273307164 &  1.273307164  & 3.2E-13  \\
			0.7 & 1.517167403 & 1.517167403 & 1.517167403 & 2.0E-12  \\ 
			0.8 & 1.776211964 & 1.776211964& 1.776211964  & 9.9E-12  \\ 
			0.9 & 2.053033451 & 2.053033451& 2.053033451 & 3.9E-11 \\ 
			1.0 & 2.350402387 &  2.350402387 & 2.350402387  & 1.3E-10  \\ 	\noalign{\smallskip}\hline
	\end{tabular}}
	\label{tab3}
\end{table}

\begin{table}[htbp]
	\centering
	\caption{Maximum absolute error for $w_{1}$, $w_{2}$ and $w_{3}$ of example 3}
	{\begin{tabular}{llll} \hline\noalign{\smallskip}
			$N$ &  $E_{1N,\infty}$ &  $E_{2N,\infty}$ &  $E_{3N,\infty}$ \\ \noalign{\smallskip}\hline\noalign{\smallskip}
			5 & 9.9E-03 & 1.2E-03 &  4.0E-04 \\ 
			10 & 3.0E-07& 2.3E-08 &  5.0E-08 \\
			15 & 8.1E-13& 7.1E-13 &  5.7E-15 \\   \noalign{\smallskip}\hline
	\end{tabular}}
	\label{tab4}
\end{table}

\begin{figure}[H]
	\centering
	\subfloat[absolute error for $w_{1}(v)$]{%
		\resizebox*{6.8cm}{!}{\includegraphics{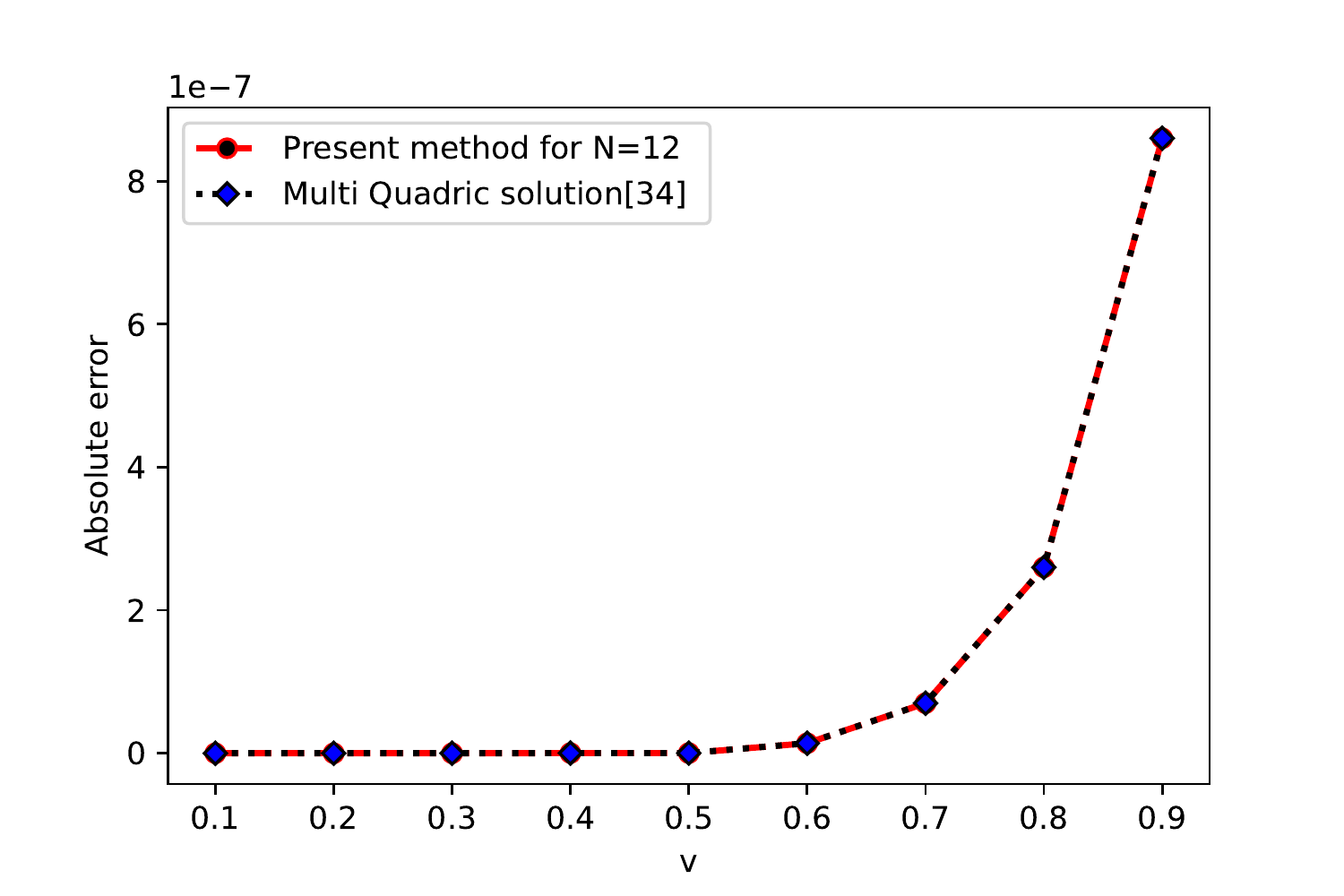}}}\hspace{3pt}
	\subfloat[absolute error for $w_{2}(v)$]{%
		\resizebox*{6.8cm}{!}{\includegraphics{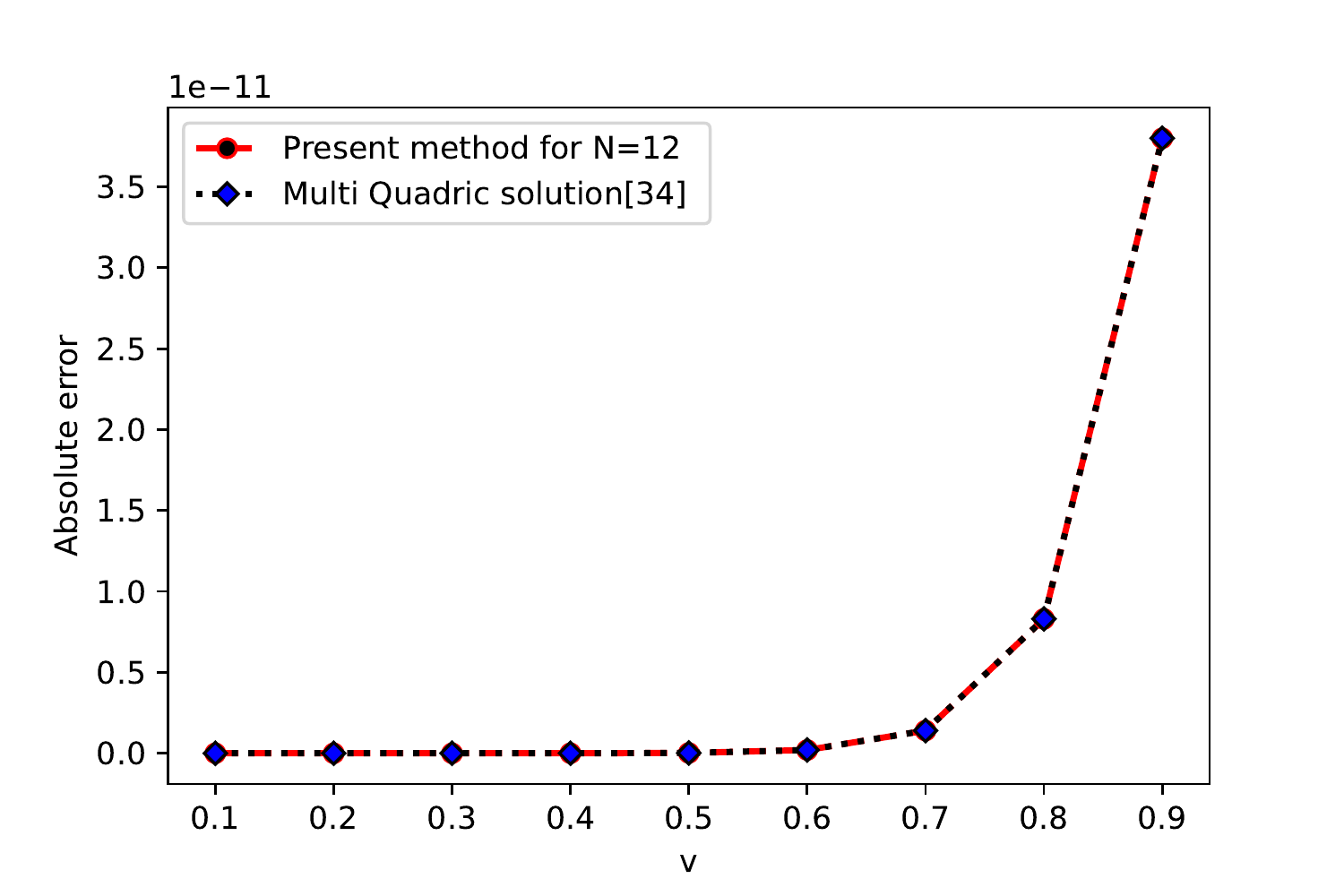}}}\hspace{3pt}
	\subfloat[absolute error for $w_{3}(v)$]{%
		\resizebox*{6.8cm}{!}{\includegraphics{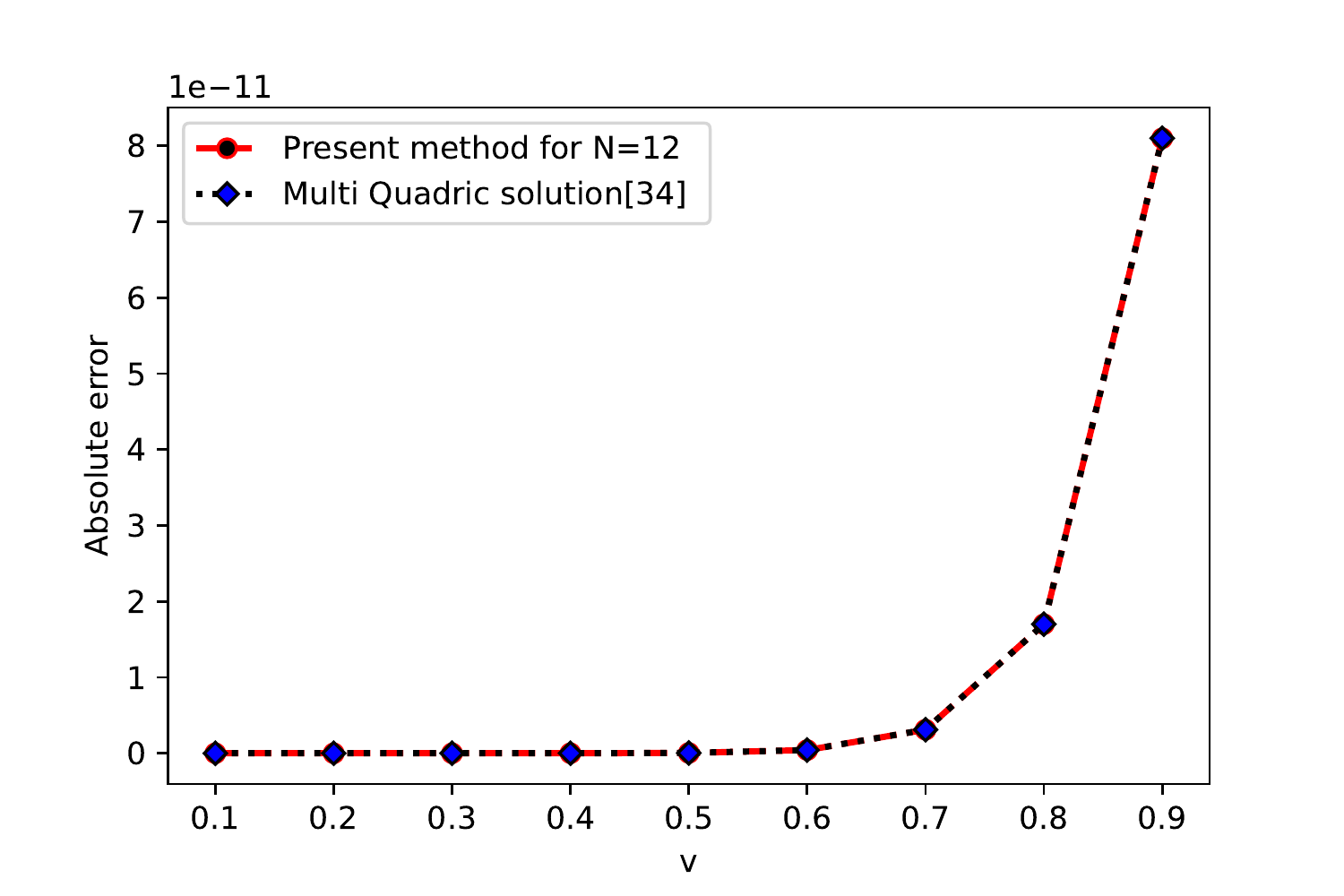}}}
	\caption{Maximum absolute error for $w_{1}$, $w_{2}$ and $w_{3}$ of example 3.} 
	\label{fig2}
\end{figure}

Tables [\eqref{tab1},\eqref{tab2},\eqref{tab3}] exhibit numerical results for $N=12$ for $w_1$, $w_2$, and $w_3$ using the proposed technique. The absolute error of numerical results obtained from the present approach is plotted in figure \eqref{fig2} against the multiquadric method \cite{VanaA11} and the maximum absolute error of numerical results obtained from the present method is shown in Table \eqref{tab4}. Furthermore, the absolute error and maximum absolute error both decline as $N$ increases, demonstrating that the method is convergent.

\subsection{Example 4}
Consider the differential algebraic equation  that studies the position of particle on a circular track, which is termed as mechanical control problem \cite{Liu03,Sand02}
\begin{eqnarray} \label{33}
	\nonumber \frac{d^{2}w_{1}}{dv^{2}}&=&2w_{2}+w_{1} w_{3}, \\
	\nonumber \frac{d^{2}w_{2}}{dv^{2}}&=&-2w_{1}+w_{2} w_{3}, \\
	0&=&w^{2}_{1} +w^{2}_{2}-1,\hspace{1mm} 0 \le v \leq 1,
\end{eqnarray} 

\noindent initial conditions
\begin{equation} \label{34}
	w_{1} \left(0\right)=0,\hspace{1mm}w_{2} \left(0\right)=1,\hspace{1mm}w_{3} \left(0\right)=0 .       
\end{equation}
The exact solution is given by
\begin{equation} \label{35} 
	w\left(v\right)=\left(\begin{array}{c} {\sin{v}^{2}} \\ {\cos{v}^{2}} \\ {-4 {v}^{2}} \end{array}\right).
\end{equation} 

\begin{table}[htbp]
	\centering
	\caption{Comparison of numerical solution of $w_{1}(v)$ with exact solution and Liu \cite{Liu03} solution  for example 4 $\left(N=20\right)$ }
	{\begin{tabular}{lllll} \hline\noalign{\smallskip}
			$v$ & $w_{1}\left(v\right)$ &$w_{1,N}\left(v\right)$  &  $w_{1,N}\left(v\right)$\cite{Liu03}  &  $R_{1,N}(v)$ \\ 	\noalign{\smallskip}\hline\noalign{\smallskip}
			0.1 & 0.009999833 & 0.009999833 & 0.009999833  & 0  \\ 
			0.2 & 0.039989334 & 0.039989334 & 0.039989334  & 0  \\ 
			0.3 & 0.089878549 & 0.089878549 & 0.089878549  & 1.5E-16  \\ 
			0.4 & 0.159318207 & 0.159318207 & 0.159318207  & 1.7E-16 \\ 
			0.5 & 0.247403959 & 0.247403959 & 0.247403959  & 2.4E-14  \\ 
			0.6 & 0.352274233 & 0.352274233 & 0.352274233  & 9.3E-13  \\
			0.7 & 0.470625888 & 0.470625888 & 0.470625888  & 2.0E-11  \\ 
			0.8 & 0.597195441 & 0.597195441 & 0.597195441  & 3.0E-10  \\ 
			0.9 & 0.724287174 & 0.724287174 & 0.724287174  & 3.3E-09 \\ 
			1.0 & 0.841470985 & 0.841470984 & 0.841470984  & 2.9E-08  \\ 	\noalign{\smallskip}\hline
	\end{tabular}}
	\label{tab5}
\end{table}

\begin{table}[htbp]
	\centering
	\caption{Comparison of numerical solution of $w_{2}(v)$ with exact solution and Liu \cite{Liu03} solution  for example 4 $\left(N=20\right)$}
	{\begin{tabular}{lllll} \hline\noalign{\smallskip}
			$v$ & $w_{2}\left(v\right)$ &$w_{2,N}\left(v\right)$  &  $w_{2,N}\left(v\right)$\cite{Liu03}  &  $R_{2,N}(v)$ \\ \noalign{\smallskip}\hline\noalign{\smallskip}
			0.0 & 1.000000000 & 1.000000000 & 1.000000000  & 0  \\ 
			0.1 & 0.999950000 & 0.999950000 & 0.999950000  & 0  \\ 
			0.2 & 0.999200107 & 0.999200107 & 0.999200107  & 1.1E-16  \\ 
			0.3 & 0.995952733 & 0.995952733 & 0.995952733  & 1.1E-16  \\ 
			0.4 & 0.987227283 & 0.987227283 & 0.987227283  & 1.1E-16 \\ 
			0.5 & 0.968912422 & 0.968912422 & 0.968912422  & 1.1E-16  \\ 
			0.6 & 0.935896824 & 0.935896824 & 0.935896824  & 1.0E-14  \\
			0.7 & 0.882332859 & 0.882332859 & 0.882332859  & 4.5E-13  \\ 
			0.8 & 0.802095755 & 0.802095755 & 0.802095755  & 1.2E-11  \\ 
			0.9 & 0.689498433 & 0.689498433 & 0.689498433 & 2.4E-10 \\ 
			1.0 & 0.540302306 & 0.540302306 & 0.540302306 & 3.8E-09  \\ \noalign{\smallskip}\hline
	\end{tabular}}
	\label{tab6}
\end{table}

\begin{table}[!htb]
	\centering
	\caption{Maximum absolute error for $w_{1}$, $w_{2}$ and $w_{3}$ of example 4}
	{\begin{tabular}{llll} \hline\noalign{\smallskip}
			$N$ &  $E_{1N,\infty}$ & $E_{2N,\infty}$ & $E_{3N,\infty}$ \\ 	\noalign{\smallskip}\hline\noalign{\smallskip}
			10 & 4.4E-05 & 3.8E-04 &  0 \\ 
			15 & 2.7E-06 & 2.7E-07 &  0 \\
			20 & 2.4E-08 & 2.0E-09 &  0 \\   \noalign{\smallskip}\hline
	\end{tabular}}
	\label{tab7}
\end{table}

\begin{figure}[H]
	\centering
	\subfloat[absolute error for $w_{1}(v)$]{%
		\resizebox*{6.8cm}{!}{\includegraphics{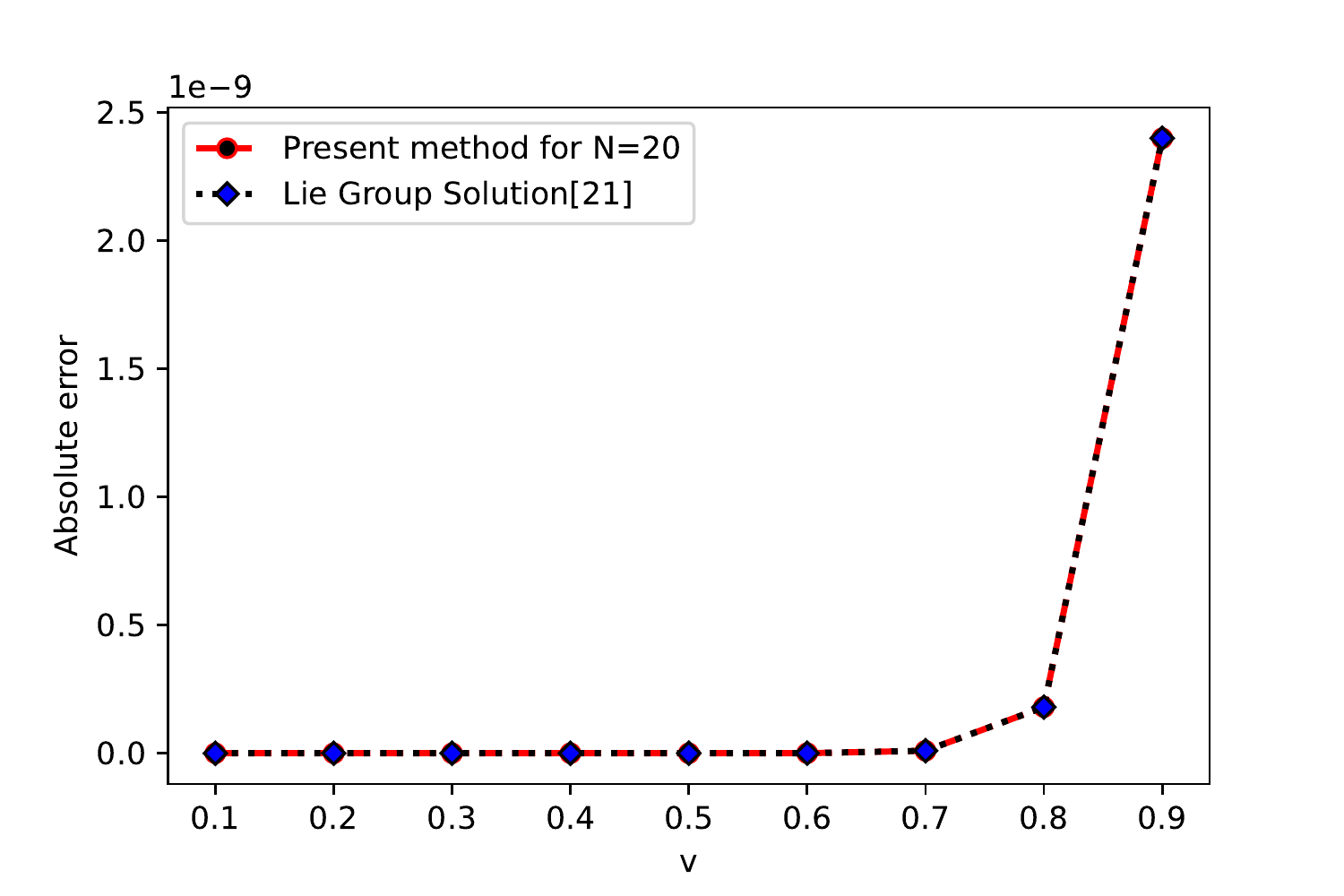}}}\hspace{3pt}
	\subfloat[absolute error for $w_{2}(v)$]{%
		\resizebox*{6.8cm}{!}{\includegraphics{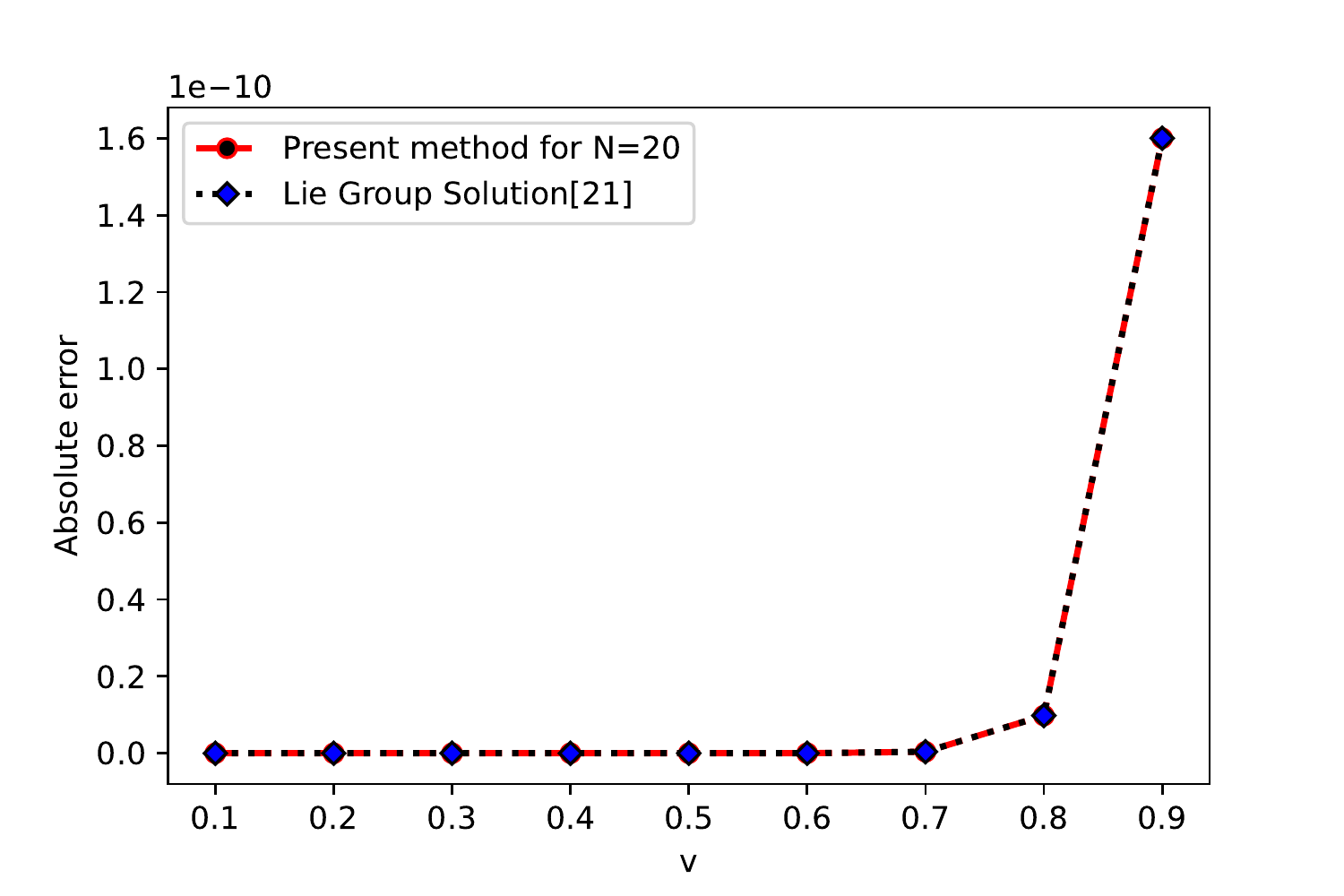}}}
	\caption{The absolute error for $w_{1}$ and $w_{2}$ of example 4.} 
	\label{fig3}
\end{figure}

Tables \eqref{tab5} and \eqref{tab6} show numerical results for $N = 20$ using the present study. In comparison to the Lie group method employed by \cite{Liu03}, the proposed method achieves a better approximation solution, as shown in Tables \eqref{tab5}-\eqref{tab6}.
Figure \eqref{fig3} depicts the comparison of absolute error results obtained using the present method and the Lie group method discussed in \cite{Liu03} and Table \eqref{tab7} lists the maximum absolute error using the present method. Table \eqref{tab7} shows the convergence of the results. As $N$ increases (see also Tables \eqref{tab5}-\eqref{tab6}), the error decreases rapidly.

\subsection{Example 5}
Consider the following nonlinear fractional differential algebraic equation
\begin{eqnarray}
	\nonumber	D^{\alpha} w_{1} + \sqrt{w_{1}}  &=& w_{2} + 2e^{2v},\\
	w_{1} - w_{2}^{2} &=& 0,
\end{eqnarray}
with initial conditions 

\begin{equation}
	w_{1}\left(0 \right) = w_{2}\left(0 \right) = 1.
\end{equation}

Exact solution for $\alpha=1$ is given by
\begin{equation}
	w_{1}\left(v \right)=e^{2v},w_{2}\left(v\right) = e^{v}.
\end{equation}

\begin{table}[htbp]
	\centering
	\caption{Comparison of numerical solution of $w_{1}(v)$ with exact solution for $\alpha=1$ and different values of $\alpha$  for example 5 $\left(N=10\right)$}\label{tab8}  
		\begin{tabular}{p{0.1\textwidth}p{0.19\textwidth}p{0.16\textwidth}p{0.16\textwidth}p{0.16\textwidth}p{0.16\textwidth}} \hline\noalign{\smallskip}
			$v$ & Exact solution for $\alpha=1$ &$w_{1,N}\left(v\right)$ for $\alpha$=1  &  $w_{1,N}\left(v\right)$ for $\alpha=0.9$  &  $w_{1,N}(v)$ for $\alpha=0.8$&  $w_{1,N}(v)$ for $\alpha=0.7$ \\ \noalign{\smallskip}\hline\noalign{\smallskip}
			0.1 &1.221402758  &1.221402758  & 1.232592481  & 1.242818926 &1.251760574 \\ 
			0.2 &1.491824698  &1.491824698  & 1.522000050  & 1.550840020 &1.577528022\\ 
			0.3 &1.822118800  &1.822118800  & 1.881442677  & 1.939997843 &1.996283505\\ 
			0.4 &2.225540928  &2.225540926  & 2.327166133  & 2.429996427 &2.531647704\\ 
			0.5 &2.718281828  &2.718281801  & 2.879128803  & 3.045186405 &3.212987402 \\ 
			0.6 &3.320116923  &3.320116716  & 3.561843540  & 3.815644105 &4.076781955 \\
			0.7 &4.055199967  &4.055198820  & 4.405409180  & 4.778498484 &5.168309651 \\ 
			0.8 &4.953032424  &4.953027348  & 5.446773813  & 5.979561723 &6.543727135\\ 
			0.9 &6.049647464  &6.049628566  & 6.731280574  & 7.475331043 &8.272630664\\ 
			1.0 &7.389056099  &7.388994709  & 8.314556977  & 9.335443038 &10.44120618 \\ \noalign{\smallskip}\hline
		\end{tabular}
\end{table}

\begin{table}[H]
	\centering
	\caption{Comparison of numerical solution of $w_{2}(v)$ with exact solution for $\alpha=1$ and different values of $\alpha$  for example 5 $\left(N=10\right)$}	\label{tab9}
		\begin{tabular}{p{0.1\textwidth}p{0.19\textwidth}p{0.16\textwidth}p{0.16\textwidth}p{0.16\textwidth}p{0.16\textwidth}} \hline\noalign{\smallskip}
			$v$ & Exact solution for $\alpha=1$ &$w_{2,N}\left(v\right)$ for $\alpha$=1  &  $w_{2,N}\left(v\right)$ for $\alpha=0.9$  &  $w_{2,N}(v)$ for $\alpha=0.8$&  $w_{2,N}(v)$ for $\alpha=0.7$ \\ \noalign{\smallskip}\hline\noalign{\smallskip}
			0.1 &1.105170918  & 1.105170918 & 1.110221816 & 1.114817889 &1.118821064\\ 
			0.2 &1.221402758  & 1.221402758 & 1.233693661 & 1.245327274 &1.255996824\\ 
			0.3 &1.349858808  & 1.349858808 & 1.371656910 & 1.392838054 &1.412898962\\ 
			0.4 &1.491824698  & 1.491824698 & 1.525505211 & 1.558844590 &1.591115024\\ 
			0.5 &1.648721271  & 1.648721271 & 1.696799638 & 1.745046381 &1.792478585 \\ 
			0.6 &1.822118800  & 1.822118800 & 1.887285142 & 1.953368506 &2.019089728 \\
			0.7 &2.013752707  & 2.013752707 & 2.098908983 & 2.185984441 &2.273321617\\ 
			0.8 &2.225540928  & 2.225540926 & 2.333842138 & 2.445345519 &2.557801126\\ 
			0.9 &2.459603111  & 2.459603103 & 2.594505315 & 2.734224962 &2.875342155\\ 
			1.0 &2.718281828  & 2.718281801 & 2.883602243 & 3.055790657 &3.228800260 \\ \noalign{\smallskip}\hline
		\end{tabular}
\end{table}

\begin{figure}[H]
	\centering{
		\subfloat[Values of $w_{1}(v)$ for different values of $\alpha$]{%
			\resizebox*{6.8cm}{!}{\includegraphics{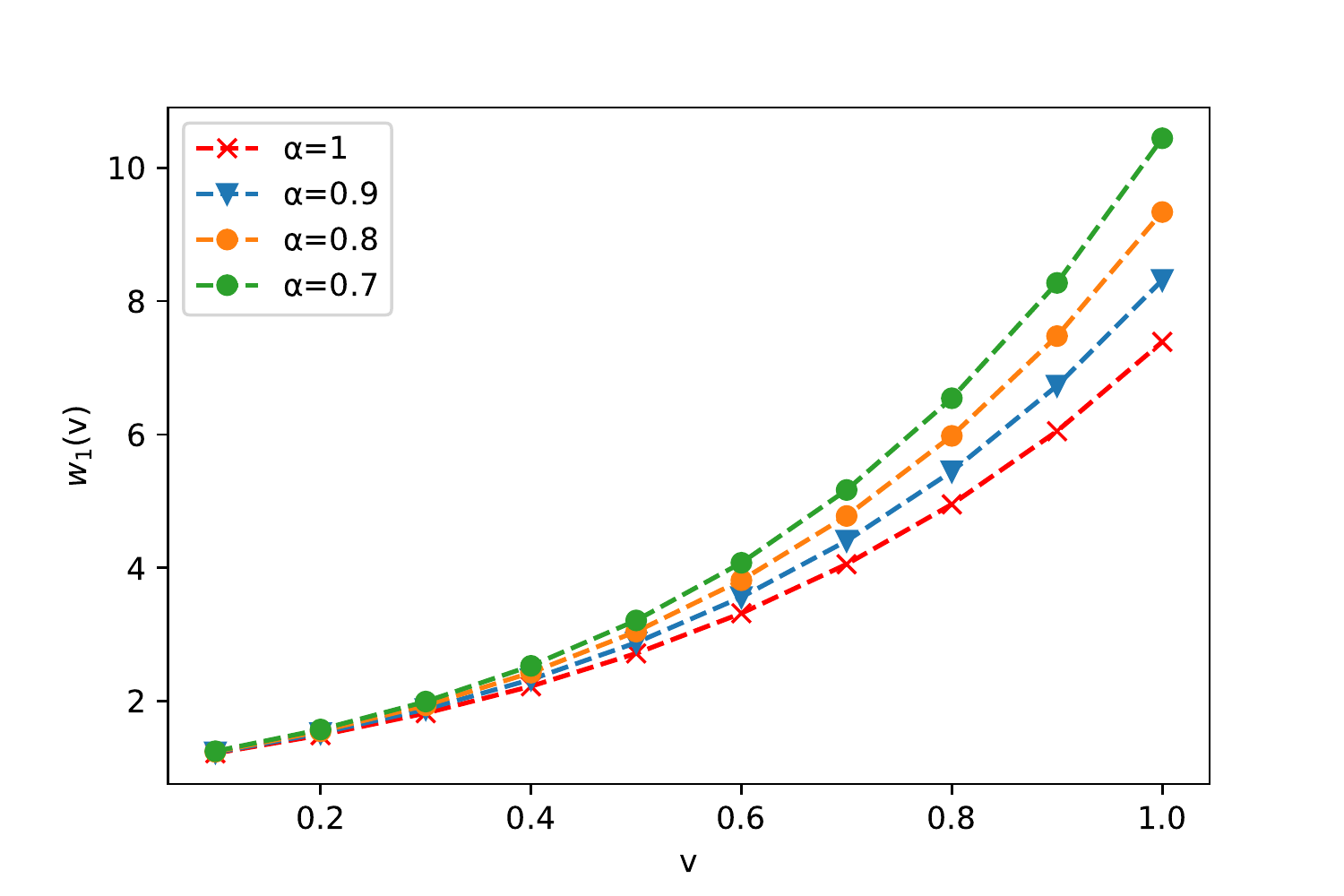}}}\hspace{3pt}
		\subfloat[Values of $w_{2}(v)$ for different values of $\alpha$]{%
			\resizebox*{6.8cm}{!}{\includegraphics{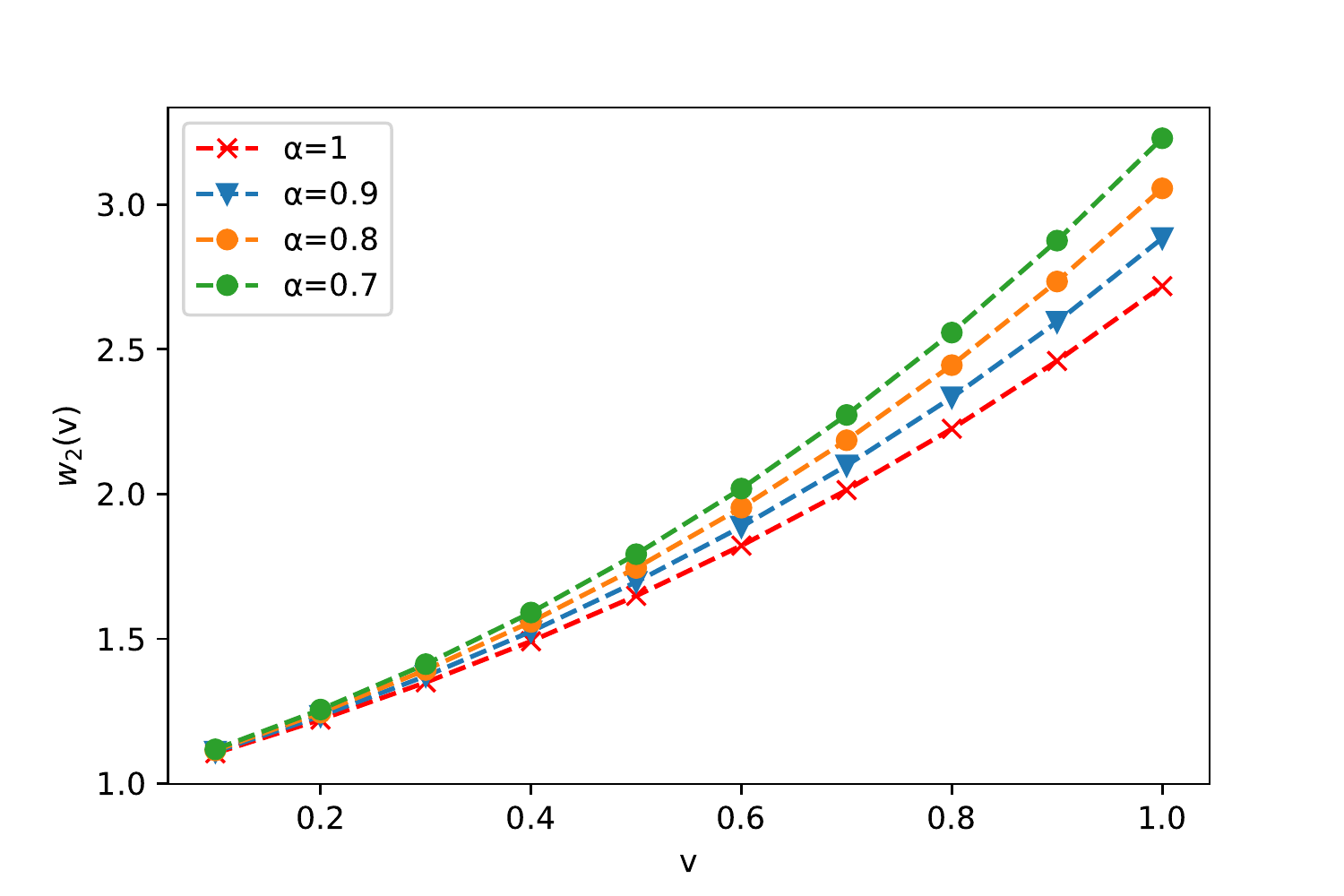}}}
		\caption{Values of $w_{1}$ and $w_{2}$ of example 5.} 
		\label{fig4}
	}
\end{figure}

Tables \eqref{tab8} and \eqref{tab9} show a comparison of numerical solutions for $w_{1}(v)$ and $w_{2}(v)$ with the exact solution and for different values of $\alpha$. Figure \eqref{fig4} shows various values for $w_{1}(v)$ and $w_{2}(v)$ for different $\alpha$. If we move away from $\alpha=1$, the solution diverges.

\section{Conclusion} \label{conclusion}
The proposed method employs the differential transform and Fa\`{a} di Bruno's formula to find solutions to nonlinear DAEs. The solutions obtained using the present technique agree quite well with the exact solution and other existing solutions. The present technique yields results that are more accurate when compared to other methods. The Fa\`{a} di Bruno's formula, which is a generalised form of the $n{th}$-order derivative of the chain rule of $f\circ g$, and the Bell polynomial are effectively used to deal with nonlinearity in the DAE problems, and gets around difficulties that previous methods had, such as the need to compute complex Adomain polynomials in ADM, the difficulty in determining the perturbation parameter, trial functions, and Lagrangian multiplier in the perturbation method, HPM, and VIM respectively, and the need to discretize variables in numerical methods that are insufficient for handling nonlinear DAEs. Other polynomials, such as Legendre polynomials, represent approximations in a given interval but not at a specific point, so only the Bell polynomial is applicable in this case. In Section \ref{applications}, the viability and applicability of the present approach are demonstrated for the mechanical control problem and the fractional DAE problem. We demonstrate that our approach is effective for a number of problems.

\smallskip

\section*{conflicts of interest}
The authors have no conflicts of interest to declare. All authors have seen and we agree with the contents of the manuscript.


\begin{thebibliography}{99}
	
	\bibitem{Aroz06} Arikoglu, A.; Ozkol, I. Solution of difference equations by using differential transform method. {\em Applied mathematics Computation}  {\bf 2006}, {\em 174}, 1216--1228. 
	
	
	\bibitem{Fatm04} Ayaz, F. Solution of the systems of differential equations by differential transform method. {\em Applied Mathematics and Computation} {\bf 2004}, {\em 147}, 547--567.	
	
	
	\bibitem{Ben15} Benhammouda, B. Solution of nonlinear higher‑index Hessenberg DAEs by Adomian polynomials and differential transform method. {\em ﻿SpringerPlus} {\bf 2015}, {\em 4}, 648--667.
	
	
	\bibitem{Ben16} Benhammouda, B. A novel technique to solve nonlinear
	higher‑index Hessenberg differential–algebraic equations by Adomian decomposition method. ﻿{\em SpringerPlus} {\bf 2016}, {\em 5}, 590--603.
	
	\bibitem{Ben22} Benhammouda, B. A New Numerical Technique for Index-3 DAEs Arising from Constrained Multibody Mechanical Systems. {\em Discrete Dynamics in Nature and Society} {\bf 2022}, {\em 2022}, 1--11.
	
	
	\bibitem{Brenan} Brenan, K.~F.; Campbell, S.~L.; Petzold, L.~R. \textit{Numerical Solution of Initial-Value Problems in Differential Algebraic Equations}, Elsevier, New York, 1989.
	
	\bibitem{Bia10} Biazar, J.; Eslami, M. Differential transform method for quadratic Riccati differential equation. {\em International Journal of Nonlinear Science} {\bf 2010}, {\em 9(4)}, 444--447.
	
	\bibitem{Carmine20} Carmine, M.~P.; Antonio, L.; Domenico, G. Stability analysis of rigid multibody mechanical systemswith holonomic and nonholonomic constraints. {\em Archive of Applied Mechanics} {\bf 2020}, {\em 90}, 1961--2005.
	
	
	\bibitem{Comt74} Comtet, L.; \textit{Advanced Combinatorics-The Art of Finite and Infinite Expansions},  D. Reidel, Dordrecht, 1974.
	
	
	\bibitem{Dehg10} Dehghan, M.; Soltanian, F.; Karbassi, S.~M. Solution of the differential algebraic equations via homotopy perturbation method and their engineering applications. {\em International Journal of Computer Mathematics} {\bf 2010}, {\em 87}, 1--10.
	
	\bibitem{Ghazwa21} Ghazwa, F.~A.; Zaboon, R.~A. Approximate solution of a reduced-type index-hessenberg differential-algebraic control system. {\em Journal of Applied Mathematics} {\bf 2021}, {\em 2021}, 1--13.
	
	\bibitem{Gha16} Ghaneai, H.; Hosseini, M.~M.; Solving differential-algebraic equations through variational iteration method with an auxiliary parameter. {\em Applied Mathematical Modelling} {\bf 2016}, {\em 40(5-6)}, 3991--4001.
	
	\bibitem{Gok12} Gokdogan, A.; Merdan, M.; Yildirim, A. The modified algorithm for the differential transform method to solution of Genesio systems.{\em Communications in Nonlinear Science and Numerical Simulation} {\bf 2012},  {\em 17(1)}, 45--51.
	
	\bibitem{Hoss06} Hosseini, M.~M. Adomian decomposition method for solution of differential-algebraic equations. {\em Journal of Computational and Applied Mathematics} {\bf 2006}, {\em 197}, 495--501.	
	
	\bibitem{JafA11} Jafari, M.~A.; Aminataei, A. Solution of differential algebraic equations via semi-analytic method. {\em Ain Shams Engineering Journal} {\bf 2011}, {\em 2}, 119--124.
	
	\bibitem{JawH18} Jawary, M.~A.; Hatif, S. A semi-analytical iterative method for solving differential algebraic equations.  {\em Ain Shams Engineering Journal} {\bf 2018},{\em 9}, 2581--2586.
	
	
	\bibitem{KaraB09} Karakoc, F.; Bereketoglu, H. Solutions of delay differential equations by using differential transform method. {\em International Journal of Computer Mathematics} {\bf 2009}, {\em 86(5)}, 914--923.
	
	
	\bibitem{KangA09} Kangalgil, F.; Ayaz, F. Solitary wave solutions for the KdV and mKdV equations by differential transform method. {\em Chaos, Solitons and Fractals} {\bf 2009}, {\em 41(1)}, 464–-472.
	
	\bibitem{Linh18} Linh, V.~H.; Truong, N.~D. Runge-Kutta methods revisited for a class of structured strangeness-free differential-algebraic equations. {\em Electronic transactions on numerical analysis ETNA} {\bf 2018}, {\em 48}, 131--155.
	
	\bibitem{Lio98} Lioen, W.~M.; Swart,J.~J. \textit{Test Set for Initial Value Problem Solvers} Centrum voor Wiskunde en Informatica(CWI), 1998.
	
	
	\bibitem{Liu03} Liu, C. Solving nonlinear differential algebraic equations by an implicit lie group method. {\em Journal of Applied Mathematics} {\bf 2013}, {\em 2013}, 1--8.
	
	\bibitem{Gir06} Methi, G. Solution of Differential Equations Using Differential Transform Method. {\em Asian Journal of Mathematics \& Statisstics} {\bf 2016}, {\em 9}, 1--5.
	
	\bibitem{Gir19} Methi, G.; Kumar, A. Numerical Solution of Linear and Higher-order Delay Differential Equations using the Coded Differential Transform Method. {\em Computer and Research Modelling} {\bf 2020}, {\em 11(6)}, 1091--1099.
	
	\bibitem{Nedia07} Nedialko,N.~S; Pryce, J.~D. Solving differential-algebraic equations
	by Taylor series. {\em Journal of Numerical Analysis,
		Industrial and Applied Mathematics}  {\bf 2007}, {\em 1(1)}, 1--30. 
	
	\bibitem{Puk82} Pukhov, G.E. Differential transforms and circuit theory. {\em Circuit theory and Applications} {\bf 1982}, {\em 10}, 265--276.
	
	\bibitem{Ravi08}RaviKanth A.S.V.; Aruna, K. Solution of singular two-point boundary value problems using differential transformation method. {\em Physics Letters A} {\bf 2008}, {\em 372}, 4671--4673.
	
	\bibitem{Reb17} Rebenda, J.; \v{S}marda, Z. A differential transformation approach for solving functional differential equations with multiple delays. {\em Communications in Nonlinear Science and Numerical Simulations} {\bf 2017}, {\em 48}, 246--257.
	
	\bibitem{Reb18} Rebenda, J. An application of Bell polynomials in numerical solving of nonlinear differential equations. {\em Aplimat Proceedings} {\bf 2018}, {\em 2018}, 1--10.
	
	\bibitem{Reb19} Rebenda, J.; \v{S}marda, Z. Numerical algorithm for nonlinear delayed differential systems of $n$th order. {\em  Advances in Difference Equations} {\bf 2019}, {\em 26}, 1--13.
	
	\bibitem{Re19} Rebenda, J. Application of Differential Transform to Multi-Term Fractional Differential Equations with Non-Commensurate Orders. {\em  Symmetry} {\bf 2019}, {\em 11 (11)}, 1390--1399.
	
	\bibitem{Reb20} Rebenda, J.; P\'{a}t\'{i}kov\'{a}, Z. Differential Transform Algorithm for Functional Differential Equations with Time-Dependent Delays. {\em Complexity} {\bf 2020}, {\em 2020}, 1--12.
	
	
	\bibitem{Sand02} Sand, J. On implicit Euler for high-order high-index DAEs. {\em Applied Numerical Mathematics}  {\bf 2002}, {\em 42(1--3)}, 411--424.
	
	
	\bibitem{Stri05} Striebel, M.; Gunther, M. A charge oriented mixed multirate method for a special class of index-1 network equations in chip design. {\em Applied Numerical Mathematics} {\bf 2005}, {\em 53(2--4)}, 489--507.
	
	
	\bibitem{VanaA11} Vanani, S.~K.; Aminataei, A. Numerical solution of differential algebraic equations using a multiquadric approximation scheme. {\em Mathematical and Computer Modelling} {\bf 2011}, {\em 53(5--6)}, 659--666.
	
	\bibitem{Wang01} Wang, H.; Song, Y. Regularization methods for solving differential-algebraic equations. {\em Applied Mathematics and Computation} {\bf 2001}, {\em 119}, 283--296.
	
	\bibitem{Yang20} Yang, L.; Kai, S. Solving power system differential algebraic equations using the differential transformation. {\em IEEE transaction on power} {\bf 2020}, {\em 35(3)}, 1--9.
	
\end{thebibliography}
\end{document}